\documentclass[pdflatex,sn-mathphys-num]{sn-jnl}
\usepackage{subfigure,graphicx}
\usepackage{graphicx}%
\usepackage{multirow}%
\usepackage{amsmath,amssymb,amsfonts}%
\usepackage{amsthm}%
\usepackage{mathrsfs}%
\usepackage[title]{appendix}%
\usepackage{xcolor}%
\usepackage{textcomp}%
\usepackage{manyfoot}%
\usepackage{booktabs}%
\usepackage{algorithm}%
\usepackage{algorithmicx}%
\usepackage{algpseudocode}%
\usepackage{listings}%

\usepackage{subfigure,graphicx}
\usepackage{graphicx}%
\usepackage{multirow}%
\usepackage{amsmath,amssymb,amsfonts}%
\usepackage{amsthm}%
\usepackage{mathrsfs}%
\usepackage[title]{appendix}%
\usepackage{xcolor}%
\usepackage{textcomp}%
\usepackage{manyfoot}%
\usepackage{booktabs}%
\usepackage{algorithm}%
\usepackage{algorithmicx}%
\usepackage{algpseudocode}%
\usepackage{listings}%
\usepackage{threeparttable}
\usepackage{diagbox}
\usepackage{supertabular}
\usepackage{supertabular}
\usepackage{booktabs}
\newcommand{\tabincell}[2]{\begin{tabular}{@{}#1@{}}#2\end{tabular}}
\usepackage{algorithm}
\usepackage{algpseudocode}
\usepackage{soul}

\usepackage{amsmath}
\usepackage{amssymb}
\usepackage{amsthm}
\usepackage{bm}
\usepackage{makecell}
\usepackage{booktabs}
\usepackage{color}
\usepackage[utf8x]{inputenc}
\usepackage{epsfig}
\usepackage{setspace}  
\usepackage{epstopdf}
\usepackage{flafter}
\usepackage{float}
\usepackage{ifthen}
\usepackage{multirow}
\usepackage{makecell} 
\usepackage{pifont}
\usepackage{natbib}
\usepackage{subfigure}
\usepackage{times}
\usepackage{threeparttable}
\usepackage{url}  
\usepackage{subfigure,graphicx}
\usepackage{slashed}
\usepackage{url}  
\usepackage{ifthen}  
\usepackage[T1]{fontenc}
\usepackage{graphicx}

\theoremstyle{thmstyleone}%
\theoremstyle{thmstyletwo}%

\theoremstyle{thmstylethree}%

\begin{document}

	\title[Article Title]{Efficient quaternion CUR method for low-rank approximation to quaternion matrix}
	
	\author[1]{\fnm{Pengling} \sur{Wu}}\email{wupengling@163.com}
		\author*[1]{\fnm{Kit } \sur{Ian Kou}}\email{kikou@um.edu.mo}
	\author[2]{\fnm{Hongmin} \sur{Cai}}\email{hmcai@scut.edu.cn}
		\author[3]{\fnm{Zhaoyuan} \sur{Yu}}\email{yuzhaoyuan@njnu.edu.cn}

	\affil*[1]{\orgdiv{Department of Mathematics}, \orgname{Faculty of Science and Technology}, \orgaddress{\street{University of Macau}, \city{ Macau 100190}, \country{China}}}		\affil[2]{\orgdiv{School of Computer Science \& Engineering}, \orgname{South China University of	Technology}, \orgaddress{\street{Guangzhou 510006}, \country{China}}}	\affil[3]{\orgdiv{Department}, \orgname{School of Geography}, \orgaddress{\street{Nanjing Normal University}, \city{ Nanjing 210023},  \country{China}}}


\abstract{	The low-rank quaternion matrix approximation has been successfully applied in many applications involving signal processing and color image processing. However, the cost of quaternion models for generating low-rank quaternion matrix approximation is sometimes considerable due to the computation of the quaternion singular value decomposition (QSVD), which limits their application to real large-scale data. To address this deficiency, an efficient quaternion matrix CUR (QMCUR) method for low-rank approximation is suggested, which provides significant acceleration in color image processing. We first explore the QMCUR approximation method, which uses actual columns and rows of the given quaternion matrix, instead of the costly QSVD. Additionally, two different sampling strategies are used to sample the above-selected columns and rows. Then, the perturbation analysis is performed on the QMCUR approximation of noisy versions of low-rank quaternion matrices. Extensive experiments on both synthetic and real data further reveal the superiority of the proposed algorithm compared with other algorithms for getting low-rank approximation, in terms of both efficiency and accuracy.}

\maketitle
\vspace{-0.5cm}
\section{Introduction}
Quaternion \cite{hamilton1866elements} as a mathematical concept was originally introduced by Hamilton in 1843. As an extension of complex numbers, a quaternion number consists of one real part and three imaginary parts. A quaternion matrix is a generalization of a complex matrix in quaternion algebra. By now, quaternions and quaternion matrices have a wide range of applications in signal processing \cite{cheng2018generalized,ell2014quaternion,LEBIHAN20041177}, machine learning \cite{6138313,ZHOU2023106234}, and particularly in color image processing \cite{ZHANG2023108971, 9066976,9994765,wang2022complex,Zhang2024complex}. By encoding the red, green, and blue channel pixel values of a color image on the three imaginary parts of quaternion matrices, this method perfectly fits the color image structure and effectively preserves the inter-relationship between the color channels \cite{qi2022quaternion}. 

As an emerging mathematical tool, low-rank quaternion matrix approximation (LRQA) has attracted much attention in the field of color image processing, such as color face recognition \cite{7468463, 8878967}, color image inpainting \cite{doi:10.1137/22M1476897,jia2019robust,li2023randomized}, and color image denoising \cite{gai2015denoising, HUANG2022108665}. A variety of LRQA-based variants are achieved by using different rank approximation regularizers \cite{9782722,chen2022color,10004995}. Notably, Chen $et\ al.$ \cite{8844978} proposed the quaternion nuclear norm (QNN)-based LRQA for color image denoising and inpainting. This method fully utilizes the high correlation among RGB channels by extending low-rank matrix approximation into the quaternion domain. Yu $et\ al.$ \cite{YU2019283} further extended the weighted nuclear norm minimization (WNNM) into the quaternion domain and proposed the quaternion-based WNNM method for color image restoration. Additionally, in \cite{YANG2021103335} and \cite{YANG202282}, the authors employed quaternion truncated nuclear norm and logarithmic norm to achieve a more accurate low-rank approximation. However, a key limitation of these methods for generating low-rank quaternion matrix approximation is that they need to compute the quaternion singular value decompositions (QSVD) in each iteration and suffer from computational deficiency, especially for large-scale data. 

In addition to rank minimization, recent studies have utilized quaternion matrix decomposition and randomized techniques to improve the performance of LRQA \cite{9623370, ren2022randomized, li2023randomized}. For instance, Miao $et\ al.$ \cite{9204671} suggested the matrix factorization for the target quaternion matrix, followed by three quaternion-based bilinear factor matrix norm factorization methods for low-rank quaternion matrix completion, which can avoid expensive calculations. Liu $et\ al.$ \cite{doi:10.1137/21M1418319} presented a randomized QSVD algorithm for low-rank matrix approximation. The randomized QSVD algorithm reduces the computational cost compared to traditional QSVD for large-scale data. Is there any method with lower computational cost for low-rank approximation of a quaternion matrix?

In recent years, the matrix CUR (MCUR) method \cite{GOREINOV19971,10.3389/fams.2018.00065,10190742} for fast low-rank approximation of real matrices has been actively investigated because of its ability to efficiently handle large-scale problems. The MCUR method approximates a low-rank matrix by directly sampling a subset of columns and rows from the original matrix and representing it as a product of three small-scale matrices. Owing to the random column/row selection strategy, the MCUR decomposition shows great potential for reducing the computational costs and preserving the properties of the original data matrix compared with other methods \cite{drineas2008relative,wang2013improving}. 

To further enhance the approximation performance and boost the computational efficiency of LRQA, we consider the efficient quaternion matrix CUR (QMCUR) method for low-rank approximation. Precisely, the QMCUR approximation of a quaternion matrix is obtained by utilizing the actual rows and columns of the original quaternion matrix. Additionally, we employ two different sampling strategies to sample the above-selected columns and rows. In summary, the introduction of the QMCUR approximation and the utilization of sampling strategies offer promising avenues for achieving low-rank approximation. The main contributions of this paper are twofold:

$\bullet$ We consider the QMCUR method for low-rank approximation of a quaternion matrix. This approach helps to reduce computational costs and improve the precision of low-rank approximation.

$\bullet$ The perturbation error bound of the proposed approximation method is studied and demonstrated under the spectral norm. Experimental results for synthetic data and color images demonstrate that the proposed QMCUR approximation method achieves a substantial speed-up while maintaining accuracy compared with other methods. 

The rest of the paper is organized as follows. The main notation and preliminary for quaternions are presented in Section 2. Section 3 presents the proposed quaternion CUR approximation method and analyzes its perturbation error bound. We provide numerical experiments to  illustrate the effectiveness of the proposed method in Section 5. Finally, some concluding remarks will be given in Section 6.
\section{Preliminaries}

\subsection{Notation}
Throughout this paper, $\mathbb{R}$ and $\mathbb{H}$ respectively denote the real and quaternion spaces. The set of all positive integers is denoted by $\mathbb{N}$, and the symbol $[n]$ represents the set of integers $\{1, \dots, n\}$ for any $n\in \mathbb{N}$. A scalar, vector, and matrix are written as $a$, $\mathbf{a}$, and $\mathbf{A}$, respectively. We use $\mathbf{A}(I,:)$ and $\mathbf{A}(:,J)$ to denote the row and column submatrices with indices $I$ and $J$, respectively. Here, $(\cdot)^*, (\cdot)^T, (\cdot)^H,$ and $(\cdot)^{\dagger}$ represent the conjugate, transpose, conjugate transpose, and Moore-Penrose inverse, respectively. $|\cdot|$ and $\|\cdot\|_F$ represent the absolute value or modulus and the Frobenius norm, respectiv

\subsection{Quaternion and quaternion matrices}
The set of quaternions $\mathbb{H}$ is a linear space over $\mathbb{R}$, with an ordered basis of 1, $\mathbf{i}$, $\mathbf{j}$, and $\mathbf{k}$. Here, $\mathbf{i}$, $\mathbf{j}$, and $\mathbf{k}$ are three imaginary units with the multiplication laws: $\mathbf{i}^2 = \mathbf{j}^2 = \mathbf{k}^2 = \mathbf{ijk} = -1$, $\mathbf{ij} = -\mathbf{ji} = \mathbf{k}$, $\mathbf{jk} = -\mathbf{kj} = \mathbf{i}$, $\mathbf{ki} = -\mathbf{ik} = \mathbf{j}$. An element \(a\) of the set of quaternions \(\mathbb{H}\) is of the form \(a = a_0 + a_1\mathbf{i} + a_2\mathbf{j} + a_3\mathbf{k}\) $\in$ \(\mathbb{H}\) with $a_0, a_1, a_2, a_3 \in \mathbb{R}$. In particular, $a$ is called a pure quaternion if the real component $a_0$ equals 0. The conjugate and modulus of \(a\) are defined as \(a^* = a_0 - a_1\mathbf{i} - a_2\mathbf{j} - a_3\mathbf{k}\) and \(|a| = \sqrt{a_0^2 + a_1^2 + a_2^2 + a_3^2}\), respectively.

For an arbitrary quaternion matrix $\mathbf{A}= \mathbf{A}_0 + \mathbf{A}_1\mathbf{i} + \mathbf{A}_2\mathbf{j} + \mathbf{A}_3\mathbf{k}= (a_{st})\in\mathbb{H}^{m\times n}$ with $\mathbf{A}_0, \mathbf{A}_1, \mathbf{A}_2, \mathbf{A}_3\in\mathbb{R}^{m\times n}$, $\mathbf{A}^H= (a^*_{ts})\in\mathbb{H}^{n\times m}$ denotes the conjugate transpose of $\mathbf{A}$. Moreover, the Frobenius norm of a quaternion matrix $\mathbf{A}\in\mathbb{H}^{m\times n}$ is $\|\mathbf{A}\|_F = \sqrt{\text{Tr}(\mathbf{A}^H\mathbf{A})}= \sqrt{\sum_{s=1}^m\sum_{t=1}^n|a_{st}|^2}$, where $\text{Tr}(\cdot)$ is the trace operator, and the spectral norm $\|\mathbf{A}\|_2$ is defined as $\|\mathbf{A}\|_2=\max_{\mathbf{x}\in\mathbb{H}^n,\|\mathbf{x}\|_2 = 1}\|\mathbf{A}\mathbf{x}\|_2$.

The quaternion singular value decomposition (QSVD) \cite{ZHANG199721} of a quaternion matrix $\mathbf{A}\in\mathbb{H}^{m\times n}$ is defined as $\mathbf{A}=\mathbf{W}\mathbf{\Sigma}\mathbf{V}^{H}$, where $\mathbf{W}\in\mathbb{H}^{m\times m}$ and $\mathbf{V}\in\mathbb{H}^{n\times n}$ are unitary, and $\mathbf{\Sigma}=\text{diag}(\sigma_1,\sigma_2, \cdots,\sigma_l)\in\mathbb{R}^{m\times n}$ with $\sigma_i\geq 0$ denoting the $i$-th largest singular value of $\mathbf{A}$ and $l = \text{min}\{m,n\}$. The truncated QSVD of a quaternion matrix $\mathbf{A}$ with rank $k$ will be denoted by $\mathbf{A}_k=\mathbf{W}_k\mathbf{\Sigma}_k\mathbf{V}_k^{H}$, where $\mathbf{W}\in\mathbb{H}^{m\times k}$ and $\mathbf{V}\in\mathbb{H}^{n\times k}$ are the $m \times k$ and $n \times k$ submatrices of $\mathbf{W}$ and $\mathbf{V}$, respectively, corresponding to choosing the first $k$ columns of each, and $\mathbf{\Sigma}_k$ is a $k \times k$ matrix containing the largest $k$ singular values. More information about quaternion matrices can be found in \cite{qi2022quaternion, ZHANG199721,ling2022singular}.

\section{Proposed approach}
In this section, we first present the QMCUR method, which is designed to efficiently compute a low-rank approximation of a large-scale quaternion matrix with low computational costs and comparable accuracy. Then we provide the perturbation estimates for the QMCUR approximation of the noisy version of low-rank matrices.
\begin{algorithm}[!htbp]
	\caption{Main steps of QMCUR approximation}
	\label{Algorithm}
	\begin{algorithmic}[1]
		\renewcommand{\algorithmicrequire}{\textbf{Input:}}
		\Require 
		$ \mathbf{X}\in\mathbb{H}^{m\times n}$: the quaternion matrix; $k$: predefined rank; $|I|, |J|$: number of rows and columns to sample;  $\{p_j\},\{q_i\}$: sampling probability distribution, $j=1,\dots,n,$ and $i = 1,\dots, m$.
		\State  Draw sampling column indices $J\subseteq [n]$ based on the sampling probability $p_j$ and construct $\mathbf{C}\in\mathbb{H}^{m\times |J|}$;
		\State Draw sampling row indices $I\subseteq [m]$ based on the sampling probability $q_i$ and construct $\mathbf{R}\in\mathbb{H}^{|I|\times n}$;
		\State {Compute} $\mathbf{U}$: $\mathbf{U}= \mathbf{C}^{\dagger}\mathbf{X}\mathbf{R}^{\dagger}\in\mathbb{H}^{|I|\times |J|}$; 
		\renewcommand{\algorithmicensure}{\textbf{Output:}}
		\Ensure 
		$\mathbf{C},\mathbf{U},\mathbf{R}$ such that  ${\mathbf{X}} \approx \mathbf{C}\mathbf{U}\mathbf{R}.$
		
	\end{algorithmic}
\end{algorithm}
\subsection{QMCUR-based low-rank approximation method}
We now present the QMCUR approximation method. Let $\mathbf{X}\in\mathbb{H}^{m\times n}$ be a low-rank quaternion matrix with a predefined rank $k$. Consider row indices $I\subseteq [m]$ and column indices $J\subseteq [n]$ satisfying $|I|, |J|\geq k$. Denote a column submatrix $\mathbf{C} = \mathbf{X}(:, J)$ whose columns span the column space of $\mathbf{X}$ and a row submatrix $\mathbf{R} = \mathbf{X}(I,:)$ whose rows span the row space of $\mathbf{X}$. Then the QMCUR approximation of $\mathbf{X}$ is a product of the form 
\begin{equation}
	\mathbf{X} \approx \mathbf{CUR},\label{Eq1}\end{equation}
where $\mathbf{C}\in\mathbb{H}^{m\times |J|}$ with $|J|$ columns from the quaternion matrix $\mathbf{X}$ and $\mathbf{R}\in\mathbb{H}^{|I|\times n}$ with $|I|$ rows from the quaternion matrix $\mathbf{X}$. 
It's obvious that the core quaternion matrix $\mathbf{U}\in\mathbb{H}^{|J|\times |I|}$ should be computed to yield the smallest error. The optimal choice for the core quaternion matrix $\mathbf{U}$ in the least-squares sense is $\mathbf{U}=\mathbf{C}^{\dagger}\mathbf{X}\mathbf{R}^{\dagger}$ \cite{WANG2005665} because
\begin{equation}
	\mathbf{C}^{\dagger}\mathbf{X}\mathbf{R}^{\dagger} = \mathop{\arg\min}_{\mathbf{U}\in\mathbb{H}^{|J|\times |I|}}\|\mathbf{X}-\mathbf{CUR}\|_F^2.
\end{equation}
Note that one may replace $\mathbf{C}^{\dagger}$ in (\ref{Eq1}) by $\mathbf{C}^{\dagger}\mathbf{X}\mathbf{R}^{\dagger}$ to improve the approximation, resulting in the QMCUR approximation $\mathbf{C}\mathbf{C}^{\dagger}\mathbf{X}\mathbf{R}^{\dagger}\mathbf{R}$. Algorithm 1 summarizes the overall process of the QMCUR approximation.

To implement the QMCUR approximation method, we set the size of $\mathbf{C}$ to be $m \times k\log k$ and $\mathbf{R}$ to be $k\log k \times n$ for the QMCUR approximation method. The indices used to determine $\mathbf{C}$ and $\mathbf{R}$ are sampled using two different strategies \cite{doi:10.1137/S0097539704442702,doi:10.1137/110852310}. In the first strategy, the sampling probabilities $p_j^{\text{col}}$ and $q_i^{\text{row}}$ for each column $j$ and row $i$ of the quaternion matrix $\mathbf{X}$ are based on the Euclidean norm of the columns and rows, which are respectively defined as $p_j^{\text{col}}:=\frac{\|\mathbf{X}(:,j)\|_2^2}{\|\mathbf{X}\|_F^2}$ where $j = 1,2,\dots,n$ and $q_i^{\text{row}}:=\frac{\|\mathbf{X}(i,:)\|_2^2}{\|\mathbf{X}\|_F^2}$, where $i = 1,2,\dots,m$. This strategy is referred to as QMCUR\_length. In the second strategy, the sampling probabilities \( p_j^{\text{unif}} \) and \( q_i^{\text{unif}} \) are respectively defined as \( p_j^{\text{unif}} := \frac{1}{n} \) and \( q_i^{\text{unif}} := \frac{1}{m} \). This approach is referred to as QMCUR\_uniform.
\subsection{Perturbation Estimates for CUR Approximation}
In practical applications, quaternion matrix data is often perturbed, such as noise in pictures, which may cause huge errors. In this section, we present the perturbation analysis suggested by the QMCUR approximation described above.  Our main task will be to consider matrices in the form of\begin{equation}
	\tilde{\mathbf{X}} = \mathbf{X} + \mathbf{E},\end{equation} where $\mathbf{X}\in\mathbb{H}^{m\times n}$ has rank $k$ with $k\geq\{m,n\}$, and $\mathbf{E}\in\mathbb{H}^{m\times n}$ is an arbitrary noise quaternion matrix drawn from a certain distribution. 

Before analyzing the perturbation, we will first introduce some notation. Let $\tilde{\mathbf{C}} = \tilde{\mathbf{X}}(:,J)\in\mathbb{H}^{m\times |J|}$, $\tilde{\mathbf{R}}= \tilde{\mathbf{X}}(I,:)\in\mathbb{H}^{|I|\times n}$, and $\tilde{\mathbf{U}} = \tilde{\mathbf{C}}^{\dagger}\tilde{\mathbf{X}}\tilde{\mathbf{R}}^{\dagger}\in\mathbb{H}^{|J|\times |I|}$ for selected index sets $I$, $J$. Then we can express this as
\begin{equation}
	\tilde{\mathbf{C}} = \mathbf{C}+ \mathbf{E}, \tilde{\mathbf{R}} = \mathbf{R}+ \mathbf{E}, \tilde{\mathbf{U}} = \mathbf{U}+ \mathbf{E}, \label{cur_noise}
\end{equation}
where ${\mathbf{C}} ={\mathbf{X}}(:,J)$, ${\mathbf{R}}= {\mathbf{X}}(I,:)$, and ${\mathbf{U}} ={\mathbf{C}}^{\dagger}{\mathbf{X}}{\mathbf{R}}^{\dagger}.$ For simplicity, we write  ${\mathbf{E}_J}={\mathbf{E}}(:,J)$ and ${\mathbf{E}_I}= {\mathbf{E}}(I,:)$. The following theorem provides a perturbation estimate for the QMCUR approximation.

\noindent \textbf{Theorem 1} Let $\tilde{\mathbf{X}} = \mathbf{X} + \mathbf{E}$ for a fixed but arbitrary $\mathbf{E}\in\mathbb{H}^{m\times n}$. Using the notation in (\ref{cur_noise}), then the following holds:
\begin{equation}
	\|\mathbf{X}-\tilde{\mathbf{C}}\tilde{\mathbf{C}}^{\dagger}\tilde{\mathbf{X}}\tilde{\mathbf{R}}^{\dagger}\tilde{\mathbf{R}}\|\leq \|\mathbf{E}_I\|\|{\mathbf{X}}{\mathbf{R}}^{\dagger}\| + \|\mathbf{E}_J\|\|{\mathbf{C}}^{\dagger}{\mathbf{X}}\| + 3\|\mathbf{E}\|.
\end{equation}  
Furthermore, 
\begin{equation}
	\|\mathbf{X}-\tilde{\mathbf{C}}\tilde{\mathbf{C}}^{\dagger}\tilde{\mathbf{X}}\tilde{\mathbf{R}}^{\dagger}\tilde{\mathbf{R}}\|\leq \|\mathbf{E}\|(\|\mathbf{W}_{k,I}^{\dagger}\| + \|\mathbf{V}_{k,J}^{\dagger}\| +3).
\end{equation} 
where $\mathbf{W}_{k,I} = \mathbf{W}_{k}(I,:)\in\mathbb{H}^{|I|\times k}$ and $\mathbf{V}_{k,J} = \mathbf{V}_{k}(J,:)\in\mathbb{H}^{|J|\times k}$ for index sets $I, J$.

In a word, Theorem 3.1 shows that the error estimation in the QMCUR method is controlled by the pseudoinverses of the submatrices of the orthogonal singular vectors and is linear in the norm of the noise $\mathbf{E}$. To establish Theorem 3.1, we will introduce the following lemmas.\\
\noindent \textbf{Lemma 1} Suppose that $\mathbf{X}\in\mathbb{H}^{m\times n}$ has rank $k$ with the QMCUR approximation: $\mathbf{X}\approx \mathbf{CUR}$, where $\mathbf{C} = \mathbf{X}(:,J)$ with selected column indices $J$,  $\mathbf{R} = \mathbf{X}(I,:)$ with selected row indices $I$, and $\mathbf{U} = \mathbf{C}^\dagger\mathbf{X}\mathbf{R}^\dagger$. Then we have 
\begin{equation}
	\text{rank}(\mathbf{C}) = \text{rank}(\mathbf{R}) = \text{rank}(\mathbf{X}).
\end{equation}
\begin{proof}
	Recall that $\mathbf{C} = \mathbf{X}(:,J)$.  Thus, span($\mathbf{C}$) $\subset$ span($\mathbf{X}$). On account of the fact that $\mathbf{X}=\mathbf{C}\mathbf{C}^{\dagger}\mathbf{X}\mathbf{R}^{\dagger}\mathbf{R}$, we have span($\mathbf{X}$) $\subset$ span($\mathbf{C}$). A similar argument shows that span($\mathbf{R}^H$) $\subset$ span($\mathbf{X}^H$). The conclusion follows.
\end{proof}
\noindent \textbf{Lemma 2} The following hold:
\begin{equation}
	\|\mathbf{X}-\tilde{\mathbf{C}}\tilde{\mathbf{C}}^\dagger\tilde{\mathbf{X}}\|\leq\|\mathbf{E}_J\|\|\mathbf{C}^\dagger\mathbf{X}\| + \|\mathbf{E}\|, \ 	\|\mathbf{X}-\tilde{\mathbf{X}}\tilde{\mathbf{R}}^\dagger\tilde{\mathbf{R}}\|\leq\|\mathbf{E}_I\|\|\mathbf{X}\mathbf{R}^\dagger\| + \|\mathbf{E}\|.
\end{equation}
\begin{proof}
	Note that
	\begin{equation}
		\begin{aligned}
			\|(\mathbf{I}-\tilde{\mathbf{C}}\tilde{\mathbf{C}}^\dagger){\mathbf{C}}\|& = 	\|(\mathbf{I}-\tilde{\mathbf{C}}\tilde{\mathbf{C}}^\dagger)\tilde{\mathbf{C}} -  (\mathbf{I}-\tilde{\mathbf{C}}\tilde{\mathbf{C}}^\dagger){\mathbf{E}_J}\| \\ & \leq \|(\mathbf{I}-\tilde{\mathbf{C}}\tilde{\mathbf{C}}^\dagger)\tilde{\mathbf{C}}\| +  \|(\mathbf{I}-\tilde{\mathbf{C}}\tilde{\mathbf{C}}^\dagger){\mathbf{E}_J}\| \\ & \leq \|\mathbf{E}_J\|.
		\end{aligned}
	\end{equation}
	The final inequality holds because the first norm term is 0 by identity of the Moore–Penrose pseudoinverse and $\|\mathbf{I}-\tilde{\mathbf{C}}\tilde{\mathbf{C}}^\dagger\|_2\leq 1$ as this is an orthogonal projection operator. Since $\text{rank}(\mathbf{C}) = \text{rank}(\mathbf{X})$, we have $\mathbf{X} = \mathbf{C}\mathbf{C}^\dagger\mathbf{X}$. Then 
	\begin{equation}
		\begin{aligned}
			\|\mathbf{X}-\tilde{\mathbf{C}}\tilde{\mathbf{C}}^\dagger\tilde{\mathbf{X}}\|&\leq	\|(\mathbf{I}-\tilde{\mathbf{C}}\tilde{\mathbf{C}}^\dagger){\mathbf{X}}\| + \|\mathbf{E}\|\\ & 
			= \|(\mathbf{I}-\tilde{\mathbf{C}}\tilde{\mathbf{C}}^\dagger){\mathbf{C}}\mathbf{C}^\dagger\mathbf{X}\| + \|\mathbf{E}\| \\&
			\leq\|\mathbf{E}_J\|\|\mathbf{C}^\dagger\mathbf{X}\| + \|\mathbf{E}\|, \
		\end{aligned} 
	\end{equation}
	The second inequality is derived by mimicking the above argument. The conclusion follows.
\end{proof} 
\noindent \textbf{Lemma 3} Suppose that $\mathbf{X}\in\mathbb{H}^{m\times n}$ has rank $k$ and its truncated QSVD is $\mathbf{X}={\mathbf{W}}_k
\mathbf{\Sigma}_k{\mathbf{V}}_k^{H}.$ Let $\mathbf{C}$ and $\mathbf{R}$ be submatrices of $\mathbf{X}$ (with selected row and column indexes $I$ and $J$, respectively) and $\mathbf{U} = \mathbf{C}^\dagger\mathbf{X}\mathbf{R}^\dagger$ such that $\mathbf{X} = \mathbf{CUR}.$ Then the following hold:
\begin{equation}
	\|{\mathbf{X}}{\mathbf{R}}^{\dagger}\|= \|\mathbf{W}_{k,I}^{\dagger}\|, \ \|{\mathbf{C}}^{\dagger}{\mathbf{X}}\| = \|(\mathbf{V}_{k,J}^H)^{\dagger}\|.
\end{equation}
where $\mathbf{W}_{k,I} = \mathbf{W}_{k}(I,:)\in\mathbb{H}^{|I|\times k}$ and $\mathbf{V}_{k,J} = \mathbf{V}_{k}(J,:)\in\mathbb{H}^{|J|\times k}$ for selected index sets $I, J$.
\begin{proof}
	Begin with the fact that $\mathbf{R}={\mathbf{W}}_k(I,:)
	\mathbf{\Sigma}_k{\mathbf{V}}_k^{H}={\mathbf{W}}_{k,I}
	\mathbf{\Sigma}_k{\mathbf{V}}_k^{H}.	$
	Then we have 
	\begin{equation}
		\mathbf{A}\mathbf{R}^\dagger = {\mathbf{W}}_{k}
		\mathbf{\Sigma}_k{\mathbf{V}}_k^{H}({\mathbf{W}}_{k,I}
		\mathbf{\Sigma}_k{\mathbf{V}}_k^{H})^\dagger.
	\end{equation}
	Notice that $\mathbf{W}_{k,I}$ has full column rank, $\mathbf{\Sigma}_k$ is a $k\times k$ matrix with full rank, and $\mathbf{V}_k^{H}$ has orthogonal rows (i.e., $(\mathbf{V}_k^{H})^\dagger= \mathbf{V}_k$). In this case, 	\begin{equation}(\mathbf{W}_{k,I}
		\mathbf{\Sigma}_k{\mathbf{V}}_k^{H})^\dagger = ({\mathbf{V}}_k^{H})^\dagger	\mathbf{\Sigma}_k^{-1}\mathbf{W}_{k,I}^\dagger = \mathbf{V}_k	\mathbf{\Sigma}_k^{-1}\mathbf{W}_{k,I}^\dagger. 	\end{equation} Consequently, the following holds:
	\begin{equation}
		\begin{aligned}
			\|{\mathbf{X}}{\mathbf{R}}^{\dagger}\| &= \|{\mathbf{W}}_{k}
			\mathbf{\Sigma}_k{\mathbf{V}}_k^{H} \mathbf{V}_k	\mathbf{\Sigma}_k^{-1}\mathbf{W}_{k,I}^\dagger\| \\ &=\|	\mathbf{\Sigma}_k{\mathbf{V}}_k^{H} \mathbf{V}_k	\mathbf{\Sigma}_k^{-1}\mathbf{W}_{k,I}^\dagger\|\\&=\|	\mathbf{\Sigma}_k\mathbf{\Sigma}_k^{-1}\mathbf{W}_{k,I}^\dagger\|\\&=\|\mathbf{W}_{k,I}^\dagger\|.
	\end{aligned}	\end{equation}
	The second equality above follows from the unitary invariance of the norm \cite{ling2022singular}. Similarly, we have  \begin{equation}\|{\mathbf{C}}^{\dagger}{\mathbf{X}}\| = \|(\mathbf{V}_{k,J}^H)^{\dagger}\|,\end{equation} whereupon the conclusion follows from the fact that $ \|(\mathbf{V}_{k,J}^H)^{\dagger} = (\mathbf{V}_{k,J}^\dagger)^H$ which has the same norm as $\mathbf{V}_{k,J}^\dagger$.
	The conclusion follows.
\end{proof}
\noindent$\textbf{\textit{Proof of Theorem 1}}$.  First note that
\begin{equation}
	\begin{aligned}
		\|	\mathbf{X} - \tilde{\mathbf{C}}\tilde{\mathbf{C}}^{\dagger}\tilde{\mathbf{X}}\tilde{\mathbf{R}}^{\dagger}\tilde{\mathbf{R}}\|&=\| \mathbf{X} -\tilde{\mathbf{C}}\tilde{\mathbf{C}}^{\dagger}\tilde{\mathbf{X}} + \tilde{\mathbf{C}}\tilde{\mathbf{C}}^{\dagger}\tilde{\mathbf{X}}( {\mathbf{I}}-\mathbf{R}^{\dagger}\mathbf{R})\|\\& \leq  \|{\mathbf{X}} - \tilde{\mathbf{C}}\tilde{\mathbf{C}}^{\dagger}\tilde{\mathbf{X}}\| + \|\tilde{\mathbf{C}}\tilde{\mathbf{C}}^{\dagger}\| \|\tilde{\mathbf{X}}( {\mathbf{I}}-\mathbf{R}^{\dagger}\mathbf{R})\|\\&   \leq  \|\mathbf{X} - \tilde{\mathbf{C}}\tilde{\mathbf{C}}^{\dagger}\tilde{\mathbf{X}}\| +  \|\tilde{\mathbf{X}}( {\mathbf{I}}-\mathbf{R}^{\dagger}\mathbf{R})\|.
	\end{aligned}
\end{equation}
The second inequality follows from the fact that $\|\mathbf{C}\mathbf{C}^{\dagger}\|_F = 1$. Next, using the formula $	\tilde{\mathbf{X}} = \mathbf{X} + \mathbf{E}$, we have 
$\|\tilde{\mathbf{X}}(\tilde{\mathbf{I}}-\tilde{\mathbf{R}}^{\dagger}\tilde{\mathbf{R}})\|\leq \|\mathbf{X}(\tilde{\mathbf{I}}-\tilde{\mathbf{R}}^{\dagger}\tilde{\mathbf{R}})\|+\|\mathbf{E}\|$ since  $\|\mathbf{I}-\tilde{\mathbf{R}}^\dagger\tilde{\mathbf{R}}\|_2\leq 1$. Then applying Lemma 2, we obtain 
\begin{equation}
	\begin{aligned}
		\|\mathbf{X}-\tilde{\mathbf{C}}\tilde{\mathbf{C}}^{\dagger}\tilde{\mathbf{X}}\tilde{\mathbf{R}}^{\dagger}\tilde{\mathbf{R}}\|&\leq\|\mathbf{X} - \tilde{\mathbf{C}}\tilde{\mathbf{C}}^{\dagger}\tilde{\mathbf{X}}\| + \|\mathbf{X}(\tilde{\mathbf{I}}-\tilde{\mathbf{R}}^{\dagger}\tilde{\mathbf{R}})\| + \|\mathbf{E}\|
		\\&\leq \|\mathbf{E}_I\|\|{\mathbf{X}}{\mathbf{R}}^{\dagger}\| + \|\mathbf{E}_J\|\|{\mathbf{C}}^{\dagger}{\mathbf{X}}\| + 3\|\mathbf{E}\|.
	\end{aligned}
\end{equation}  
Moreover, from Lemma 3, it follows that 
\begin{equation}
	\|\mathbf{X}-\tilde{\mathbf{C}}\tilde{\mathbf{C}}^{\dagger}\tilde{\mathbf{X}}\tilde{\mathbf{R}}^{\dagger}\tilde{\mathbf{R}}\|\leq \|\mathbf{E}\|(\|\mathbf{W}_{k,I}^{\dagger}\| + \|\mathbf{V}_{k,J}^{\dagger}\| +3).
\end{equation}  
\section{Numerical experiments}
In this section, we conduct numerical experiments to assess the accuracy and computation time of the proposed QMCUR method for computing low-rank approximation of a quaternion data matrix. Several real-world datasets are employed to demonstrate our proposed algorithm.

$\textbf{\textit{Example 4.1}}$ In this experiment, we will demonstrate error estimation in Theorem 1 using simulation data. We randomly generate a quaternion matrix $\mathbf{X}$ as a product $\mathbf{X} = {\mathbf{W}}_k\mathbf{\Sigma}_k{\mathbf{V}}_k^{H}\in\mathbb{H}^{m\times m}$ and set its rank $k$ to be 50, where $\mathbf{W}_k\in\mathbb{H}^{m\times k}, \mathbf{V}_k\in\mathbb{H}^{m\times k}$ are unitary matrices, and $\mathbf{\Sigma}_k$ is a diagonal $k \times k$ matrix with positive diagonal elements. In addition, the random noise quaternion matrix is given by $\mathbf{E}=\mathbf{E}_0+\mathbf{E}_1\mathbf{i}+\mathbf{E}_2\mathbf{j}+\mathbf{E}_3\mathbf{k}\in\mathbb{H}^{m\times n}$, where the entries of $\mathbf{E}_t(t = 0, 1, 2, 3)$ are i.i.d. Gaussian variables with mean zero and variance $\sigma$.  Then we obtain the perturbed quaternion matrix
\begin{figure*}[!htbp]
	\centering	
	\subfigure{}{\includegraphics[width=0.4\textwidth]{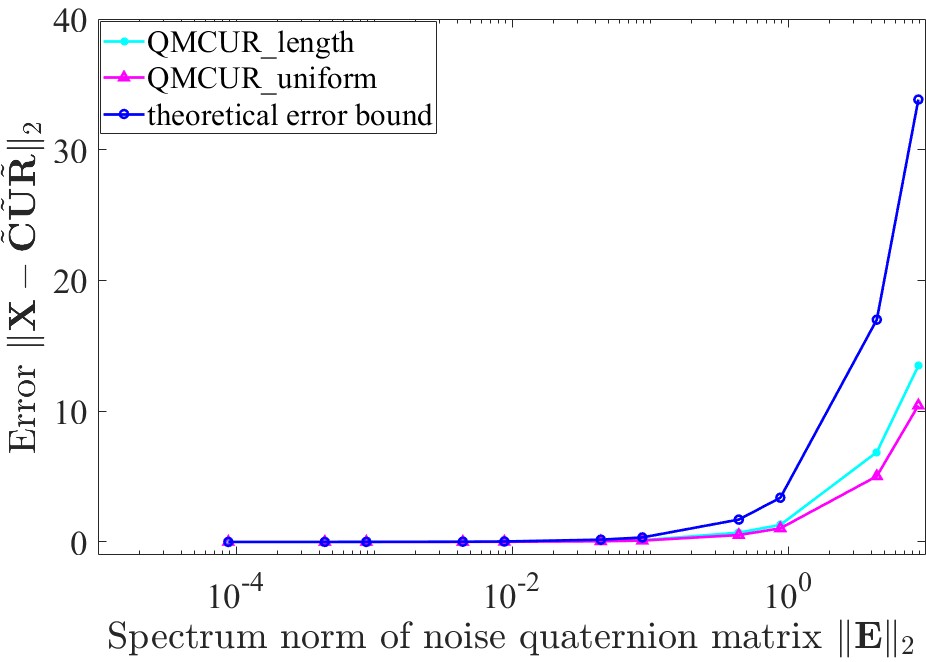}}
	\caption{$\|\mathbf{X}-\tilde{\mathbf{C}}\tilde{\mathbf{U}}\tilde{\mathbf{R}}\|_2$ vs. $\|\mathbf{E}\|_2$.} 
	\label{Figure_error}
\end{figure*}
\begin{figure*}[!htbp]
	\begin{center}	\scriptsize\setlength{\tabcolsep}{1pt}		\begin{tabular}{cccc}	
			\renewcommand\arraystretch{1}	$\sigma=10^{-1}$& 	$\sigma=10^{-4}$& 	$\sigma{-7}$ &	$\sigma =0$ \\  
\includegraphics[width=0.24\linewidth]{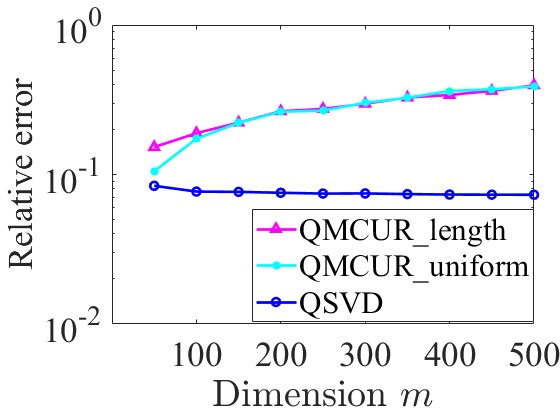}&
\includegraphics[width=0.24\linewidth]{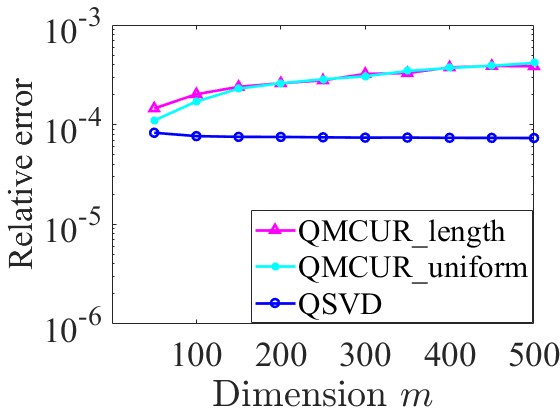}&
			\includegraphics[width=0.24\linewidth]{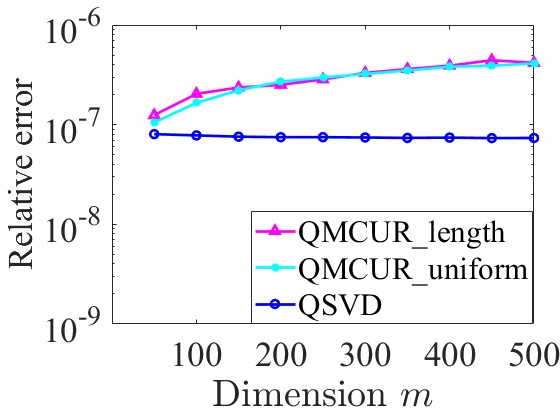}&
			\includegraphics[width=0.24\linewidth]{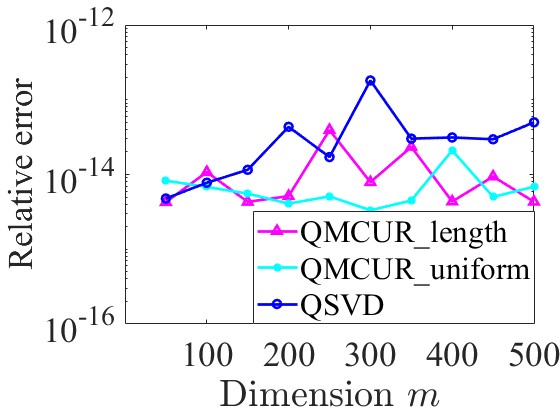}\\
			\includegraphics[width=0.24\linewidth]{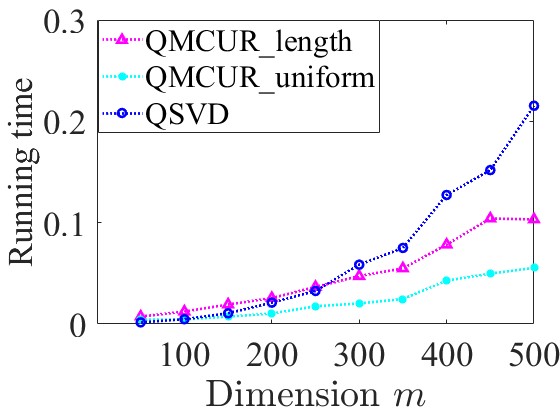}&
			\includegraphics[width=0.24\linewidth]{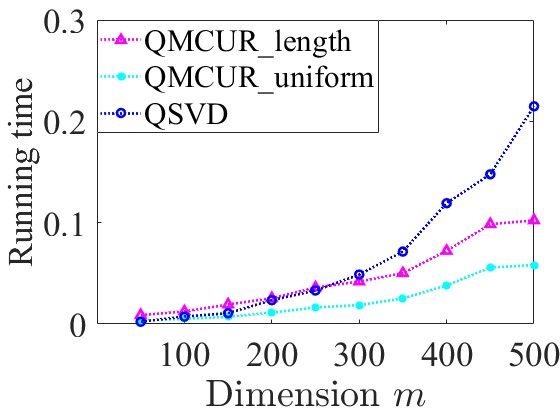}& \includegraphics[width=0.24\linewidth]{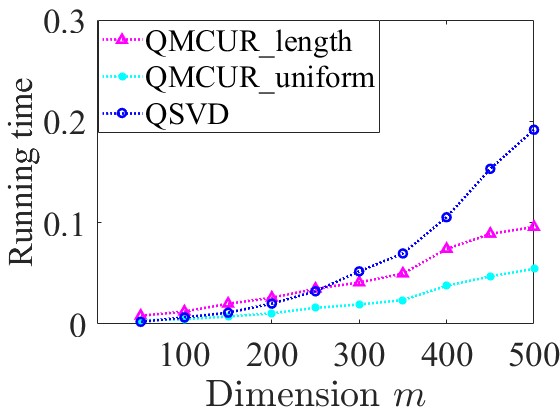}&
			\includegraphics[width=0.24\linewidth]{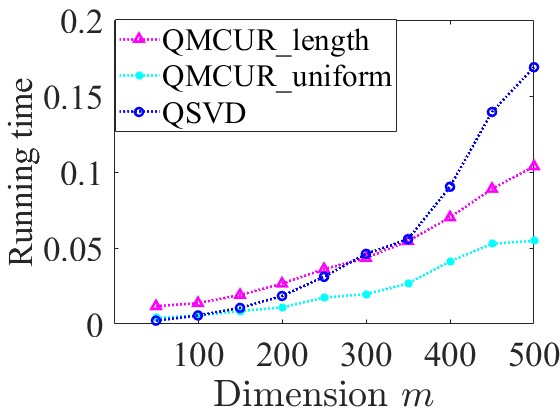} 	\end{tabular}	\caption{Comparison of low rank quaternion matrix approximation methods under different noise levels $\sigma$. Rank 10 is used in all tests and $m$ varies from 50 to 500. Top row: relative approximation errors vs. quaternion matrix dimensions. Bottom row: running time vs. quaternion matrix dimensions.}  \label{Simulat}
\end{center}\end{figure*}  
\begin{equation}
	\tilde{\mathbf{X}} = \mathbf{X}+ \mathbf{E}.
\end{equation}
\indent
In this experiment, we assume that \( m = 500 \) and \( k = 50 \). By applying Algorithm 1 to compute the rank-$k$ approximation of $\tilde{\mathbf{X}}$ with $\sigma = 10^{-1}, 10^{-2}, \dots, 10^{-6}$, Fig. \ref{Figure_error} shows the results. From this, we observe that the error $\|\mathbf{X}-\tilde{\mathbf{C}}\tilde{\mathbf{U}}\tilde{\mathbf{R}}\|_2$ increases linearly as $\|\mathbf{E}\|_2$ increases, which confirms the conclusion of Theorem 1.

$\textbf{\textit{Example 4.2}}$ In this experiment, through simulation data, we will show the efficiency and accuracy of the QMCUR method by comparison with the QSVD method. We also consider the perturbed quaternion matrix $\tilde{\mathbf{X}}$ in Example 4.1.For all methods, the relative error ${\|\mathbf{B}-\mathbf{X}\|_F}/{\|\mathbf{X}\|_F}$ and the running time (in second) are used as metrics to measure approximate solutions, where $\mathbf{B}$ is obtained by QSVD, QMCUR\_length, and QMCUR\_uniform.
\begin{figure*}[!htbp]
	\begin{center}	\scriptsize\setlength{\tabcolsep}{1pt}		\begin{tabular}{cc}		   
			\includegraphics[width=0.48\linewidth]{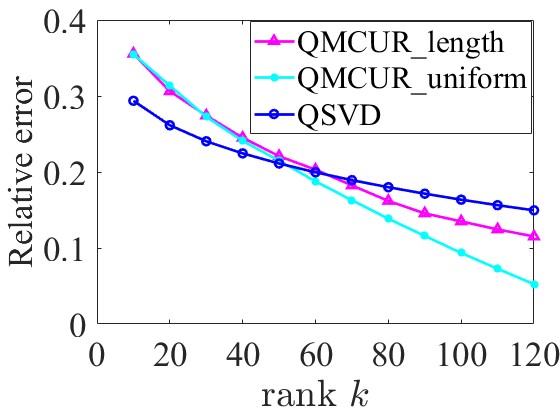}&
			\includegraphics[width=0.48\linewidth]{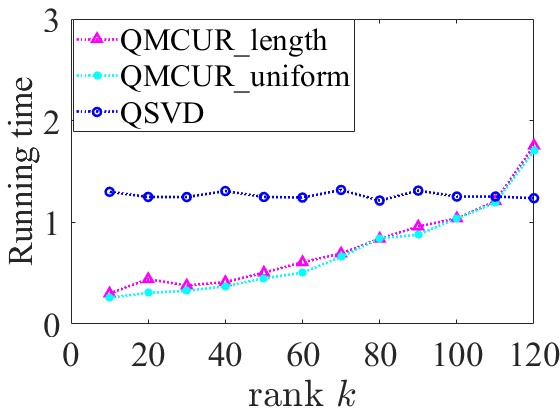}\\
			\includegraphics[width=0.48\linewidth]{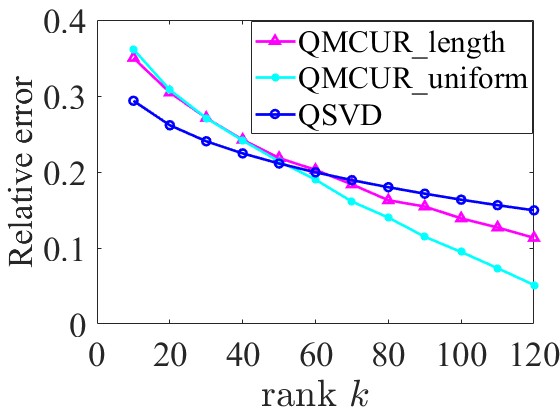}&
			\includegraphics[width=0.48\linewidth]{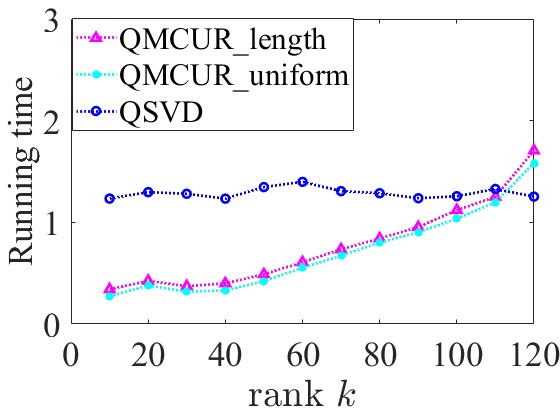}
		\end{tabular}	\caption{Numerical simulation results of QSVD and QMCUR with different sampling strategies to color images with different rank $k$.  From top to bottom: Image01 and Image02, respectively.} \label{colorlow}
\end{center}  \end{figure*} 
\\\indent In Fig. \ref{Simulat}, we present the relative error and running time for the synthetic data with varying size $m$, where the rank $k$ is set to 10. One can see that both variations of the proposed QMCUR approximation are significantly faster than QSVD when $m$ is sufficiently large. In the noiseless case (i.e., when $\sigma$ = 0), the errors for QMCUR\_uniform are smaller than those obtained by the QSVD algorithm. However, when additive noise appears, the proposed methods exhibit slightly worse but still good approximation accuracy. In summary, the proposed methods achieve a substantial speed-up while maintaining accuracy compared to QSVD, especially when $m$ is sufficiently large.
\begin{figure*}[!htbp]	\begin{center}	\scriptsize\setlength{\tabcolsep}{1pt}		\begin{tabular}{cccc}		\includegraphics[width=0.24\linewidth]{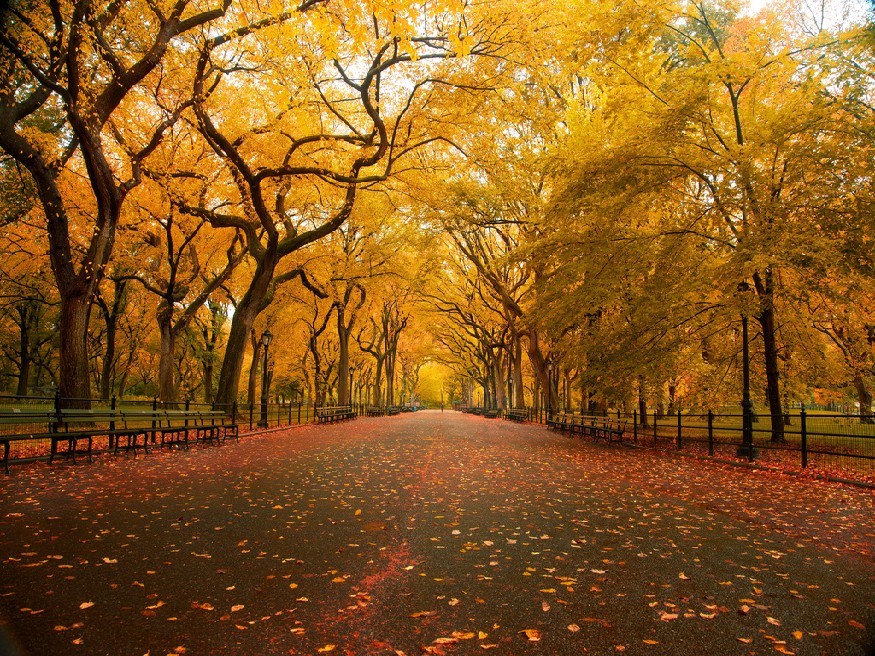}&	\includegraphics[width=0.24\linewidth]{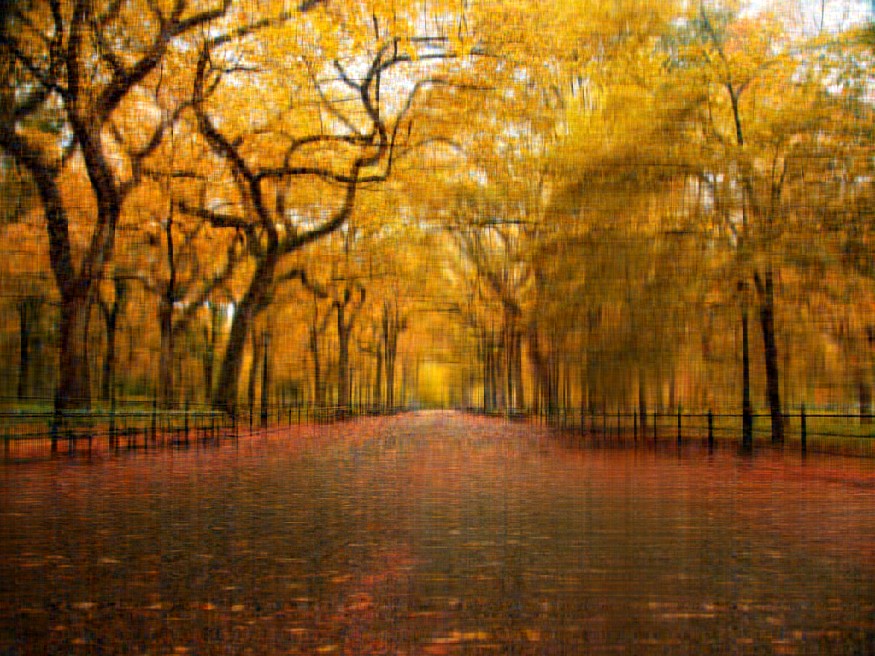}& 	\includegraphics[width=0.24\linewidth]{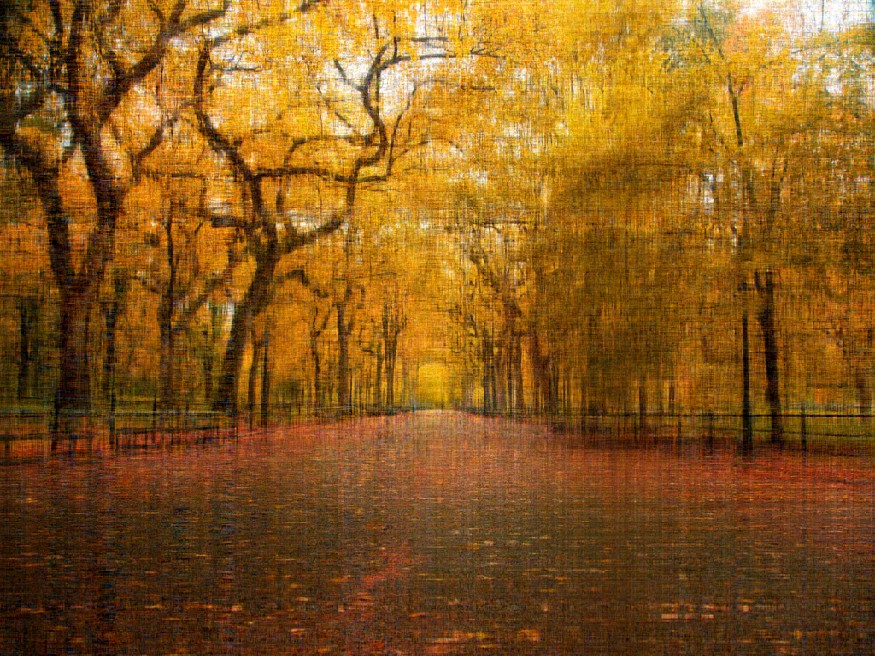}&	\includegraphics[width=0.24\linewidth]{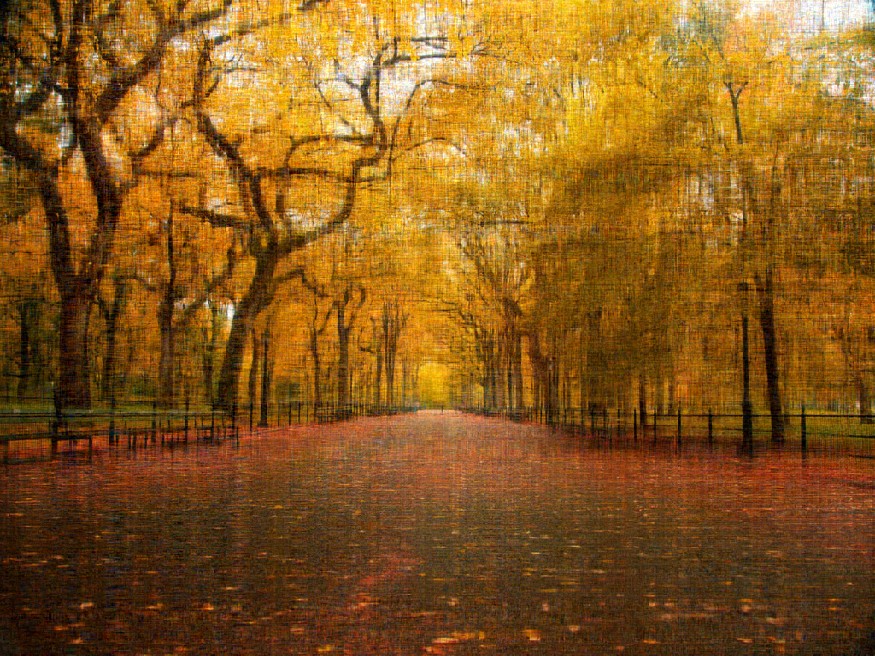}\\  	\renewcommand\arraystretch{1}  &$k=40$ & $k=40$ & $k$ =40 \\		
& 	\includegraphics[width=0.24\linewidth]{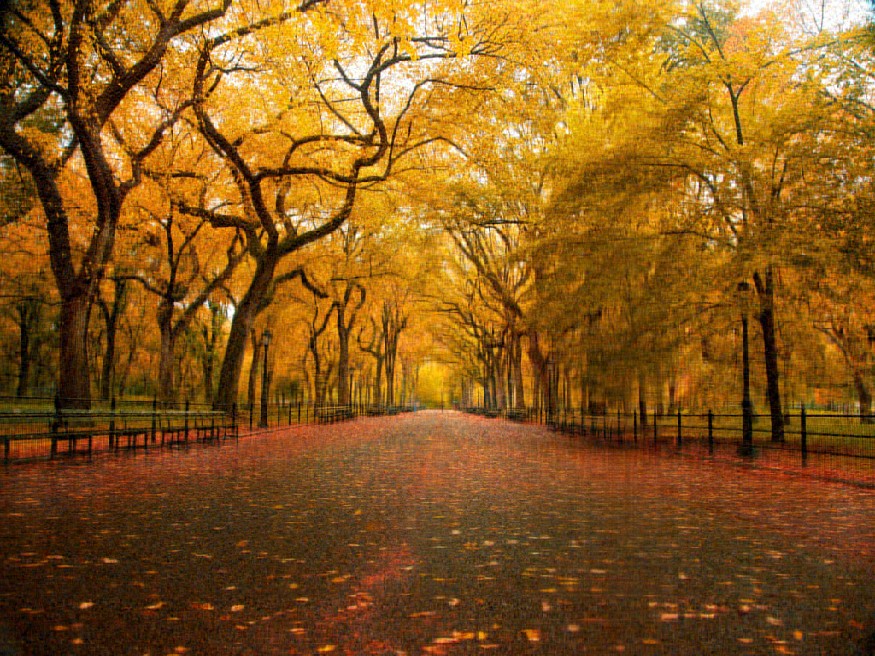}&
		\includegraphics[width=0.24\linewidth]{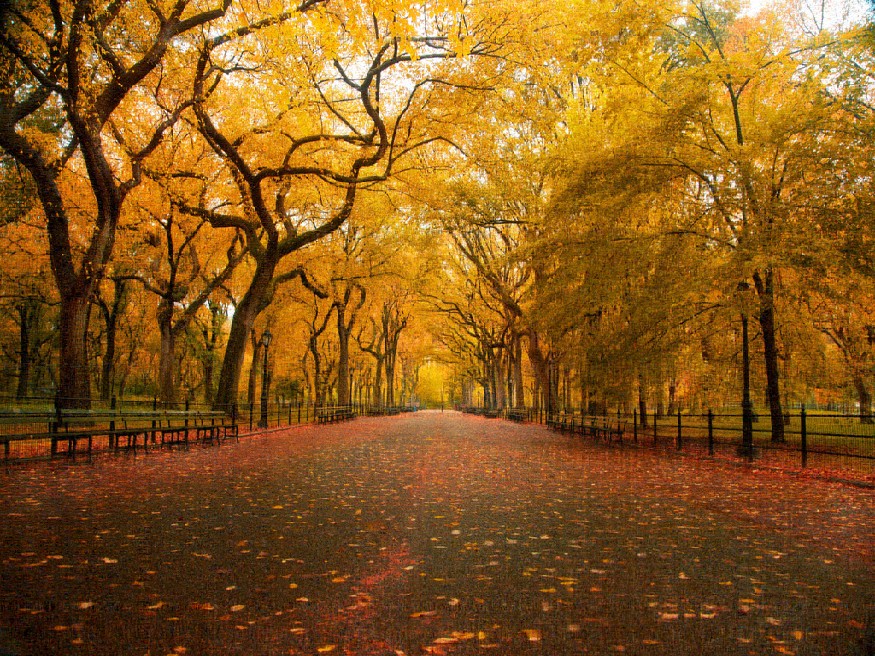}&
		\includegraphics[width=0.24\linewidth]{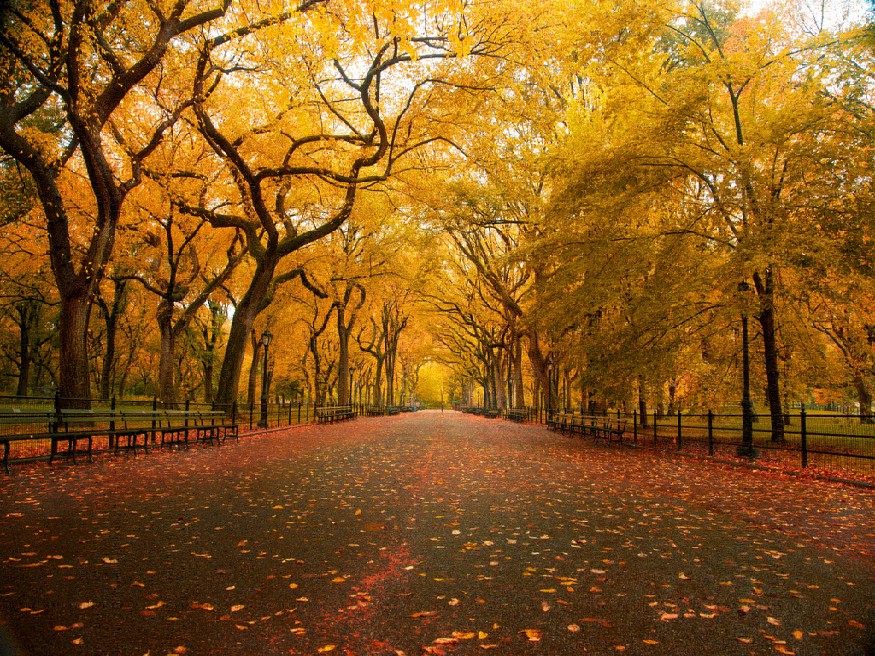}\\
		\renewcommand\arraystretch{1}  &$k=100$& $k=100$ & $k=100$ \\
	\includegraphics[width=0.24\linewidth]{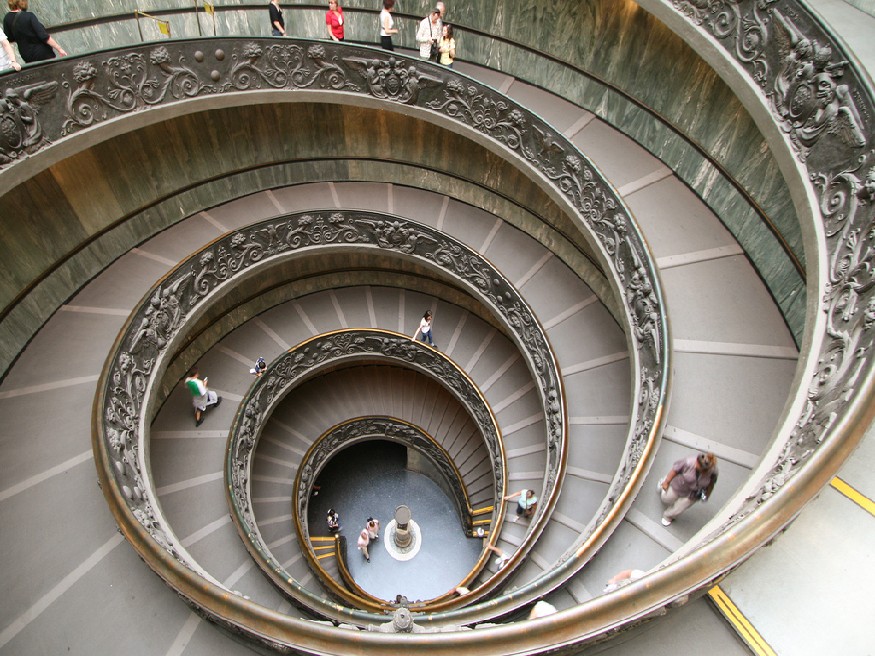}&	\includegraphics[width=0.24\linewidth]{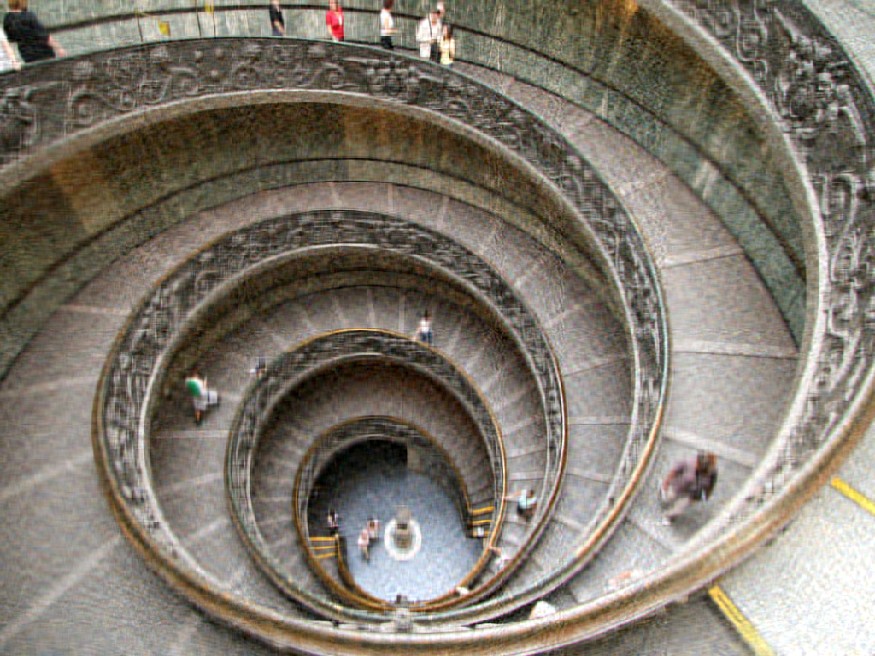}&	\includegraphics[width=0.24\linewidth]{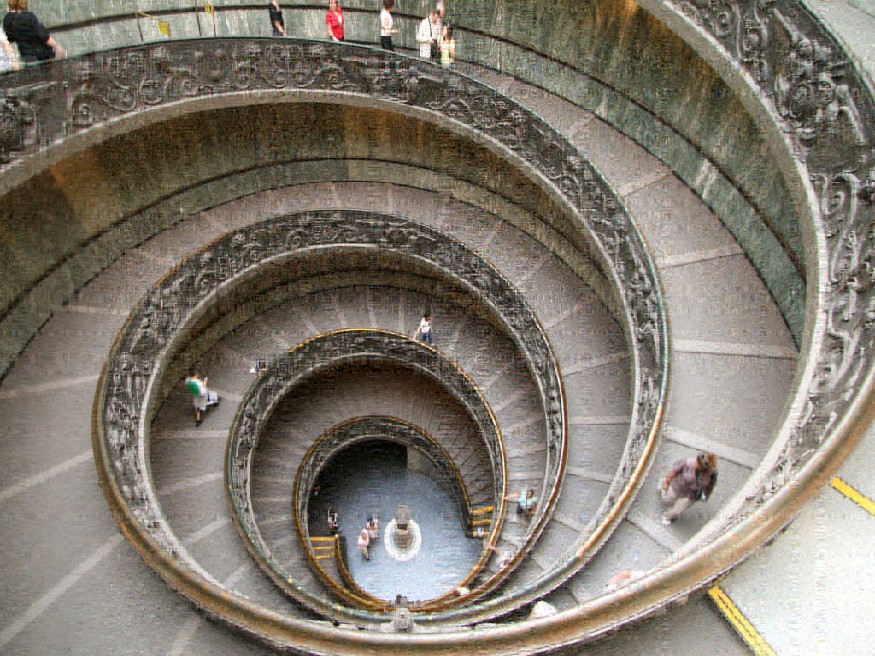}&
		\includegraphics[width=0.24\linewidth]{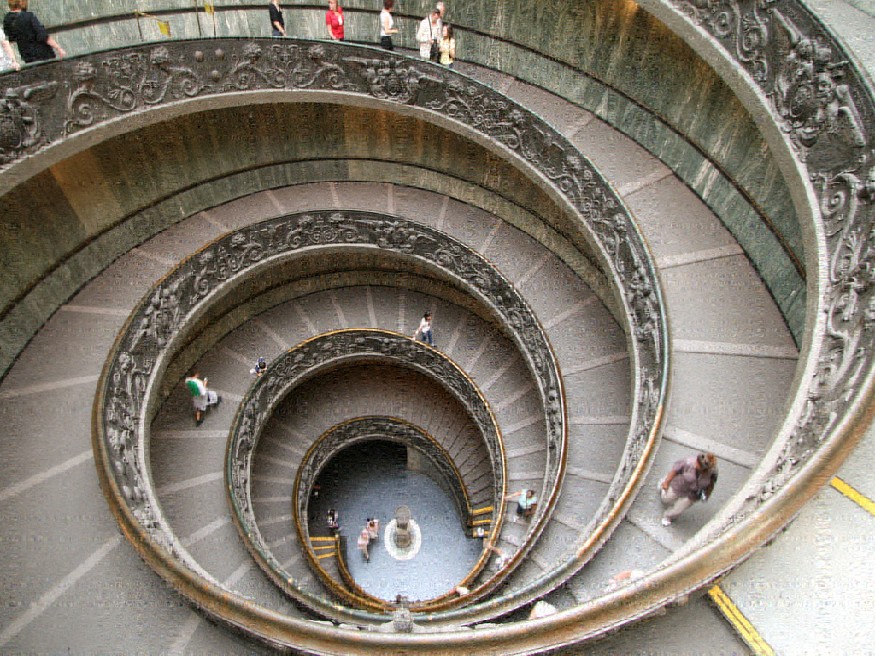} \\
		\renewcommand\arraystretch{1}   &$k=70$&$k=70$ &$k=70$\\	
	&	 \includegraphics[width=0.24\linewidth]{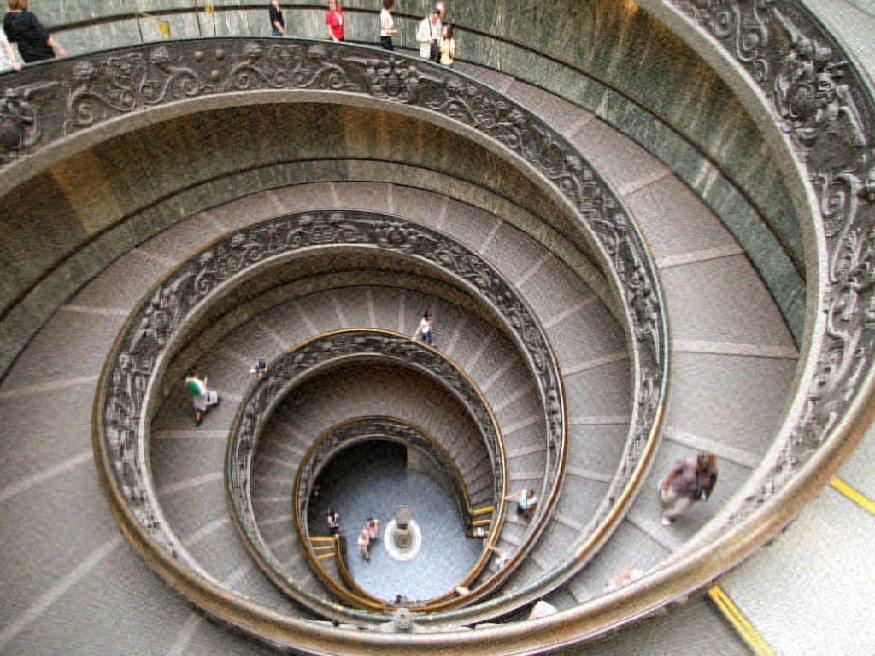}&
	\includegraphics[width=0.24\linewidth]{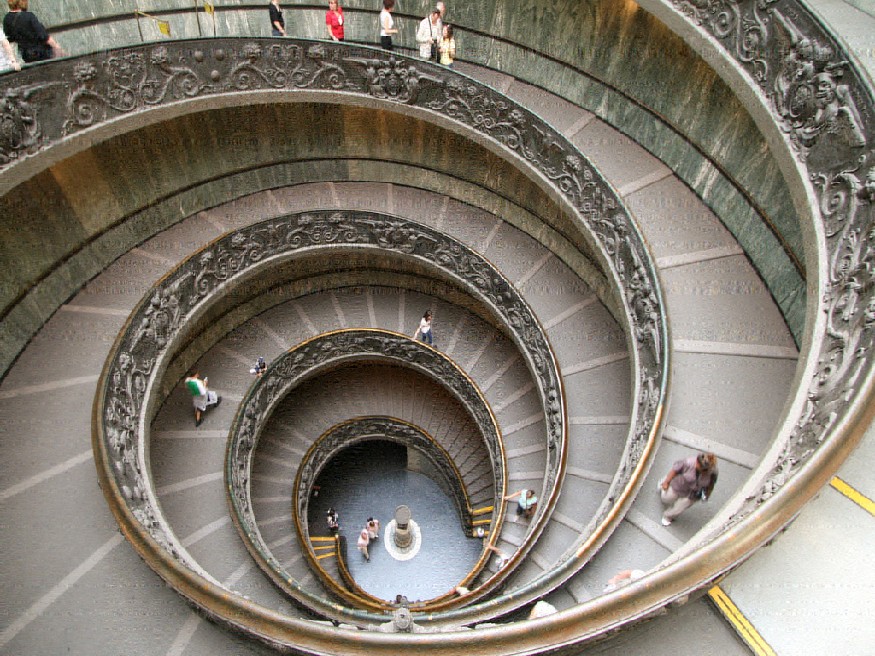}&
	\includegraphics[width=0.24\linewidth]{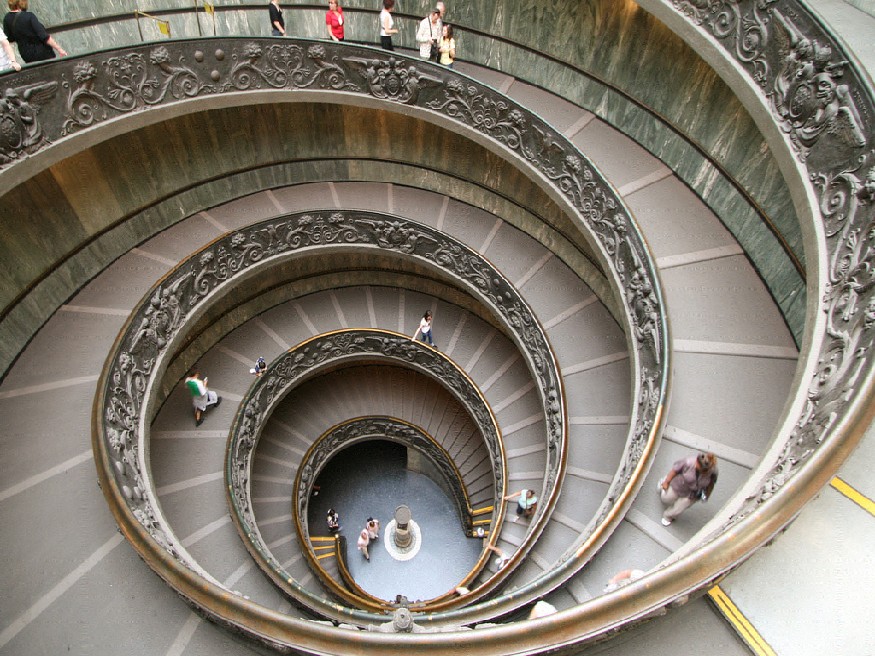}\\
		\renewcommand\arraystretch{1}   &$k=110$&$k=110$&$k=110$ \\
		\renewcommand\arraystretch{1} Original  &QSVD &QMCUR$\_$length &QMCUR$\_$uniform 
			\\\end{tabular}\caption{Low-rank color image reconstruction. From top to bottom: Image01 and Image02, respectively.} \label{coloshowww}\end{center}  \end{figure*} 

$\textbf{\textit{Example 4.3}}$ In this example, we evaluate the performance of the proposed QMCUR method for color image compression. The dataset used in this example is the Set27$\footnote{https://github.com/Huang-chao-yan/dataset.}$ and two images are shown in Fig. \ref{coloshowww}. Since the size of all images is 683 $\times$ 1024 $\times$ 3, each image can be represented by a pure quaternion matrix $\mathbf{X}= \mathbf{R}\mathbf{i} + \mathbf{G}\mathbf{j} + \mathbf{B}\mathbf{k}$, where $\mathbf{R}$, $\mathbf{G}$, and $\mathbf{B}$ represent the RGB (red, green, blue) channels of the color image. The elements of $\mathbf{X}$ are $x_{ij} = r_{ij}\mathbf{i} + g_{ij}\mathbf{j} + b_{ij}\mathbf{k}$, where $r_{ij}$, $g_{ij}$, and $b_{ij}$ denote the red, green, and blue pixel values at the position $(i, j)$ of the color image, respectively.

Fig. \ref{colorlow} illustrates the relative error and running time of all methods for various selected ranks. From Fig. \ref{colorlow}, we can observe that the performance of the proposed methods is superior to that of QSVD in terms of relative error when $k$ is set to a suitable size. In terms of running time, the proposed methods are much faster than QSVD. In Fig. \ref{coloshowww}, we present the reconstructed image quality for various target ranks. It is evident that there is little difference in the visual effect between each original image and the approximate image obtained by the proposed methods and the QSVD.

$\textbf{\textit{Example 4.4}}$ In this part, the proposed models are applied to the recovery of color images. The problem of recovering color images can be formulated as a quaternion matrix completion problem, aiming to find its missing entries through the following optimization.
\begin{equation} \min \limits_{{\mathbf{X}}} \ \text{rank}{({\mathbf{X}})},\ \ \text {s.t.}\ \ P_{\Omega}({\mathbf{X}})=P_{\Omega}({\mathbf{Y}}), \label{eq1} \end{equation}
where rank($\cdot$) denotes the rank function, ${\mathbf{X}}\in\mathbb{H}^{I_1\times I_2}$ and ${\mathbf{Y}}\in\mathbb{H}^{I_1\times I_2}$ represent the recovered and observed quaternion matrices, respectively. $\Omega$ represents the set of observed elements, and $P_{\Omega}(\mathbf{X})$ is a projection operator defined as $P_{\Omega}(\mathbf{X})_{ij} = \mathbf{X}_{ij}$ if $(i,j)\in\Omega$ and 0 otherwise. The rank minimization problem of a quaternion matrix can be expressed as the quaternion nuclear norm minimization (QNNM) \cite{8844978}. However, the algorithm that utilizes QNNM involves the QSVD of quaternion matrices, which results in high time complexity. The QNNM method can be more effectively replaced by the concept of quaternion decomposition. The quaternion matrix decomposition variant of formulation (\ref{eq1}) is presented as follows\begin{equation} \min \limits_{{\mathbf{X}}} \ \|P_{\Omega}({\mathbf{X}}) - P_{\Omega}({\mathbf{Y}})\|_F^2, \ \ \ \text {s.t.} \ \ \text{rank}{({\mathbf{X}})} = k, \label{eq2} \end{equation}
\begin{algorithm}[!htbp]
	\caption{Quaternion matrix-CUR approximation from partially observed entries}
	\label{Algorithm}
	\begin{algorithmic}[1]
		\renewcommand{\algorithmicrequire}{\textbf{Input:}}
		\Require
		An incomplete quaternion matrix  $\mathbf{Y}\in\mathbb{H}^{m\times n}$,  observation locations $\Omega$, rows and columns (i.e., ${I}$ and ${J}$) indices that define $\mathbf{R}\in\mathbb{H}^{|I|\times n}$ and $\mathbf{C}\in\mathbb{H}^{m\times |J|}$ respectively,  and maximal iteration $t_{max}$=200. 
		\renewcommand{\algorithmicensure}{\textbf{Output:}}
		\Ensure 
		The recovered quaternion matrix $\mathbf{X}^*.$	
		\State 
		$\mathbf{X}^0=\mathbf{Y}\in\mathbb{H}^{|I|\times |J|}$ is the observed quaternion matrix with missing pixels.
		\State   $\mathbf{M}^0\in\mathbb{H}^{|I|\times |J|}$ is a zero quaternion matrix.
		\State	$\mathbf{For}$ $t =0,1,2,\dots$ $\textbf{do}$ 
		\State \ \ \ \ \ \ \  $\mathbf{M}^{t+1} \leftarrow $ compute  QMCUR approximation of the data matrix $\mathbf{X}^{t}$ using Algorithm 1
		\State \ \ \ \ \ \ \ $\mathbf{X}^{t+1} = (\mathbf{M}^{t+1})_{\Omega^c}+ \mathbf{Y}_{\Omega}$, where $\Omega^c$ denotes the complementary set of $\Omega$
		\State \ \ \ \ \ \ \  $\mathbf{If}$ $\|\dot{\mathbf{X}}^{t+1} - \dot{\mathbf{X}}^t\|_{F}/\|\dot{\mathbf{X}}^t\|_{F}\leq 10^{-4},$ or $t\ge t_{max}$ $\mathbf{then}$\\
		\ \ \ \ \ \ \ \ \ \ \ \ \ $\mathbf{X}^* = \mathbf{X}^{t+1}$ and $\textbf{break}$ \\
		\ \ \ \ \ \ \  $\textbf{end}$\\
		$\textbf{end}$
	\end{algorithmic}
\end{algorithm}
where we assume that the unknown quaternion matrix $\mathbf{X}$ has a low-rank matrix representation. By introducing an auxiliary tensor variable $\mathbf{M}$, the optimization problem (\ref{eq2}) can be expressed as follows
\begin{equation} \min \limits_{{\mathbf{X}}} \ \| \mathbf{X} - \mathbf{M}\|_F^2, \ \ \ \text {s.t.} \ \ \text{rank}{({\mathbf{X}})} = k, \ \ P_{\Omega}({\mathbf{X}}) = P_{\Omega}({\mathbf{Y}}). \label{eq3} \end{equation}
Therefore, we can alternatively solve the optimization problem (\ref{eq2}) over the variables $\mathbf{X}$ and $\mathbf{M}$. The solution to the minimization problem (\ref{eq2}) can be approximated using the following iterative procedures.
\begin{equation}
	\mathbf{M}^t\leftarrow \mathcal{L}(\mathbf{X}^t),\label{EqX} \end{equation}
\begin{equation} \mathbf{X}^{t+1}\leftarrow \Omega\circledast\mathbf{Y} + (\mathbf{1}-\Omega) \circledast\mathbf{M}^t, \label{EqOut}	
\end{equation}
where $t$ represents the iteration number, $\mathcal{L}$ is an operator used to calculate a low-rank quaternion matrix approximation of the quaternion matrix $\mathbf{X}$, which can be achieved by the QMCUR technique, and $\mathbf{1}$ is a matrix in which all components are equal to one. The algorithm consists of two main steps: low-rank quaternion matrix approximation (Eq. \ref{EqX}) and masking computation (Eq. \ref{EqOut}). We summarize the proposed method in Algorithm 2. It starts with the initial incomplete data $\mathbf{X}^0$ and iteratively enhances the approximate solution until a stopping criterion is satisfied or the maximum number of iterations is reached.

We test our algorithm using four-color images of size 512$\times$768$\times$3 from the Kodak PhotoCD Dataset$\footnote{http://r0k.us/graphics/kodak/}$ under the random missing scenario. The missing ratio (MR) is defined as the ratio of the number of missing elements to the total number of elements. The quality of the restoration results is measured using peak signal-to-noise rate (PSNR), structural similarity (SSIM), and running time per iteration (in seconds). Several low-rank matrix completion methods are compared in the experiments, including LRQA-1 \cite{8844978}, QLNF \cite{YANG202282}, and Q-DFN \cite{9204671}. The numerical results of the color images recovered using various methods for different MRs are presented in Table \ref{tablecolor}, and the visualized results are depicted in Fig. \ref{tablecolorshow}.

\begin{table}[h]
	\renewcommand{\arraystretch}{1.2}
	\begin{center}	\setlength{\tabcolsep}{1.6mm}{	
		\caption{PSNR, SSIM, and running time  per iteration (in second) of results by different methods with different MRs on color images. The best and second best values are respectively highlighted in boldface and underlined.}\label{tablecolor}%
		\begin{tabular}{@{}c|c|ccc|ccc|ccc@{}} \hline
			\multirow{2}{*}{Data} & \multirow{2}{*}{Method} & \multicolumn{3}{c}{MR=90\%}                      & \multicolumn{3}{c}{MR=80\%}                      & \multicolumn{3}{c}{MR=70\%}                      \\
			&    & PSNR           & SSIM           & Time           & PSNR           & SSIM           & Time           & PSNR           & SSIM           & Time      	\\	\hline
			\multirow{6}[2]{*}{\tabincell{c}{ \includegraphics[width=0.08\linewidth]{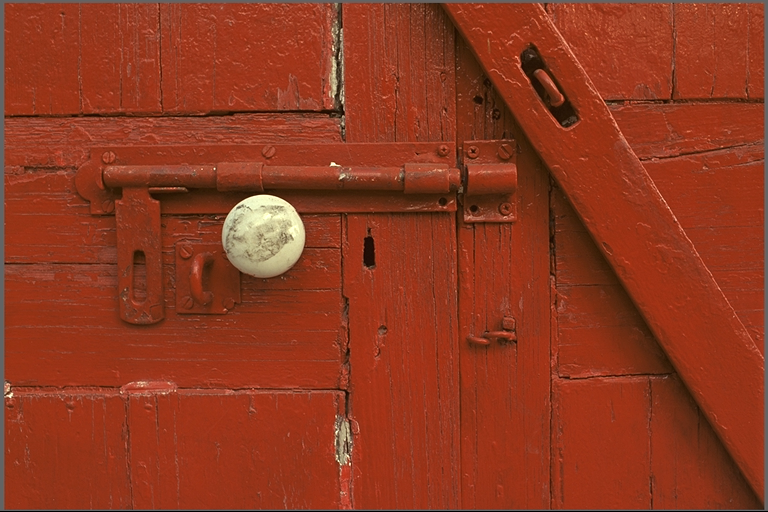} \\ {Image03  } } }	  & Observed                & 9.18           & 0.038          & 0.000          & 9.69           & 0.104          & 0.000          & 10.26          & 0.184          & 0.000          \\
			& LRQA-1                  & 24.71          & \textbf{0.954} & 0.460          & 26.73          & \underline{ 0.967}    & 0.446          & 28.07          & \underline{ 0.975}    & 0.503          \\
			& DFN                     & 24.89          & \underline{ 0.953}    & 0.348          & 27.04          & \underline{ 0.967}    & 0.280          & 28.51          & \textbf{0.976} & 0.350          \\
			& QLNF                    & 24.63          & 0.949          & 0.738          & 26.53          & 0.964          & 1.178          & 27.98          & 0.972          & 3.999          \\
			& QMCUR\_length           & \textbf{25.12} & 0.951          & \underline{ 0.226}    & \underline{ 27.20}    & \underline{ 0.967}    & \underline{ 0.257}    & \underline{ 28.63}    & \textbf{0.976} & \underline{ 0.316}    \\
			& QMCUR\_uniform          & \underline{ 24.99}    & 0.950          & \textbf{0.159} & \textbf{27.23} & \textbf{0.968} & \textbf{0.166} & \textbf{28.66} & \textbf{0.976} & \textbf{0.265}   \\	\hline		
			\multirow{6}[2]{*}{\tabincell{c}{ \includegraphics[width=0.08\linewidth]{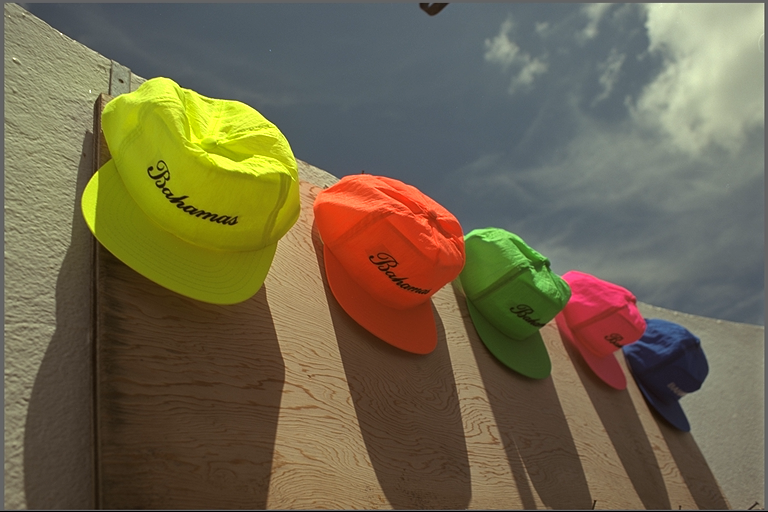} \\ {Image04  } } } & Observed                & 8.00           & 0.024          & 0.000          & 8.50           & 0.050          & 0.000          & 9.09           & 0.082          & 0.000          \\
			& LRQA-1                  & 22.79          & \underline{ 0.771}    & 0.456          & 25.50          & \textbf{0.868} & 0.498          & 27.24          & \textbf{0.905} & 0.478          \\
			& DFN                     & 23.10          & 0.757          & 0.270          & 25.96          & 0.852          & 0.350          & 27.79          & \underline{ 0.889}    & 0.326          \\
			& QLNF                    & 23.25          & 0.772          & 1.188          & 25.87          & 0.854          & 2.727          & 27.15          & 0.856          & 3.768          \\
			& QMCUR\_length           & \underline{ 23.49}    & \textbf{0.776} & \underline{ 0.184}    & \underline{ 26.23}    & \underline{ 0.858}    & \underline{ 0.321}    & \underline{ 27.88}    & \underline{ 0.889}    & \underline{ 0.346}    \\
			& QMCUR\_uniform          & \textbf{23.67} & \textbf{0.776} & \textbf{0.135} & \textbf{26.36} & \underline{ 0.848}    & \textbf{0.277} & \textbf{27.98} & 0.883          & \textbf{0.322}  \\	\hline		
			\multirow{6}[2]{*}{\tabincell{c}{ \includegraphics[width=0.08\linewidth]{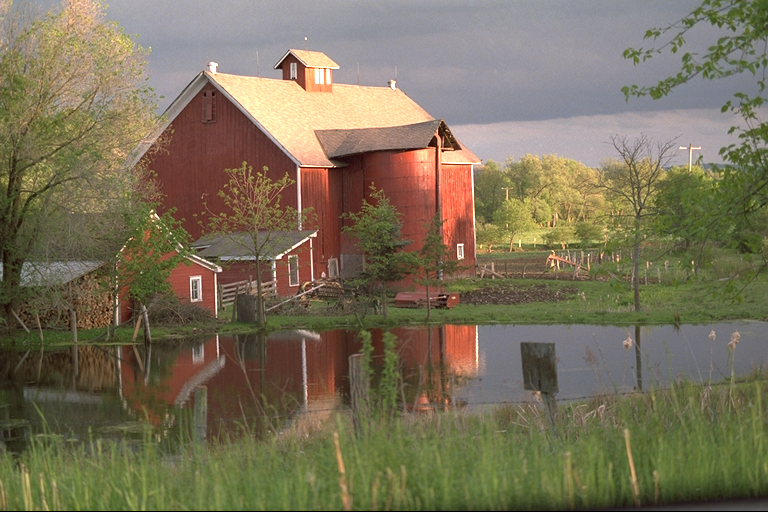} \\ {Image05 } } }  & Observed                & 7.22           & 0.022          & 0.000          & 7.73           & 0.045          & 0.000          & 8.30           & 0.072          & 0.000          \\
			& LRQA-1                  & 21.85          & \textbf{0.725} & 0.478          & 23.98          & \textbf{0.798} & 0.451          & 25.41          & \textbf{0.841} & 0.542          \\
			& DFN                     & 21.86          & \underline{ 0.702}    & 0.362          & 24.13          & \underline{ 0.783}    & 0.335          & \underline{ 25.60}    & \underline{ 0.825}    & 0.349          \\
			& QLNF                    & 21.85          & 0.690          & 1.483          & 23.83          & 0.764          & 1.583          & 18.75          & 0.648          & 2.867          \\
			& QMCUR\_length           & \underline{ 22.04}    & 0.690          & \underline{ 0.220}    & \underline{ 24.20}    & 0.779          & \underline{ 0.349}    & 25.50          & 0.824          & \underline{ 0.351}    \\
			& QMCUR\_uniform          & \textbf{22.15} & 0.688          & \textbf{0.165} & \textbf{24.22} & 0.769          & \textbf{0.283} & \textbf{25.62} & 0.822          & \textbf{0.247}       \\	\hline		
			\multirow{6}[2]{*}{\tabincell{c}{ \includegraphics[width=0.08\linewidth]{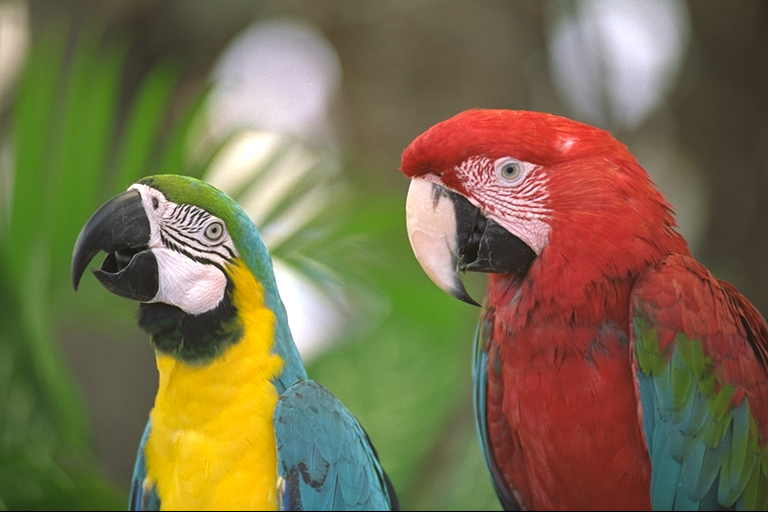} \\ {Image06  } } }  & Observed                & 7.20           & 0.026          & 0.000          & 7.71           & 0.058          & 0.000          & 8.29           & 0.095          & 0.000          \\
			& LRQA-1                  & 21.75          & 0.815          & 0.434          & 24.85          & \underline{ 0.894}    & 0.512          & 26.84          & \underline{ 0.927}    & 0.501          \\
			& DFN                     & 22.58          & \textbf{0.830} & 0.286          & 25.46          & \textbf{0.900} & 0.318          & \textbf{27.57} & \textbf{0.930} & 0.362          \\
			& QLNF                    & 22.37          & 0.810          & 1.310          & 25.02          & 0.883          & 2.789          & 26.44          & 0.907          & 4.110          \\
			& QMCUR\_length           & \underline{ 22.76}    & 0.819          & \underline{ 0.230}    & \underline{ 25.50}    & 0.892          & \underline{ 0.283}    & 27.51          & 0.923          & \underline{ 0.341}    \\
			& QMCUR\_uniform          & \textbf{22.98} & \underline{ 0.829}    & \textbf{0.168} & \textbf{25.60} & \underline{ 0.894}    & \textbf{0.227} & \underline{ 27.55}    & 0.924          & \textbf{0.327} \\	\hline	
	\end{tabular}}	\end{center}  
\end{table}
\begin{figure*}[!htbp]	\vspace{-0.3cm}	\begin{center}	\scriptsize\setlength{\tabcolsep}{1pt}		\begin{tabular}{ccccccc}		   \includegraphics[width=0.13\textwidth]{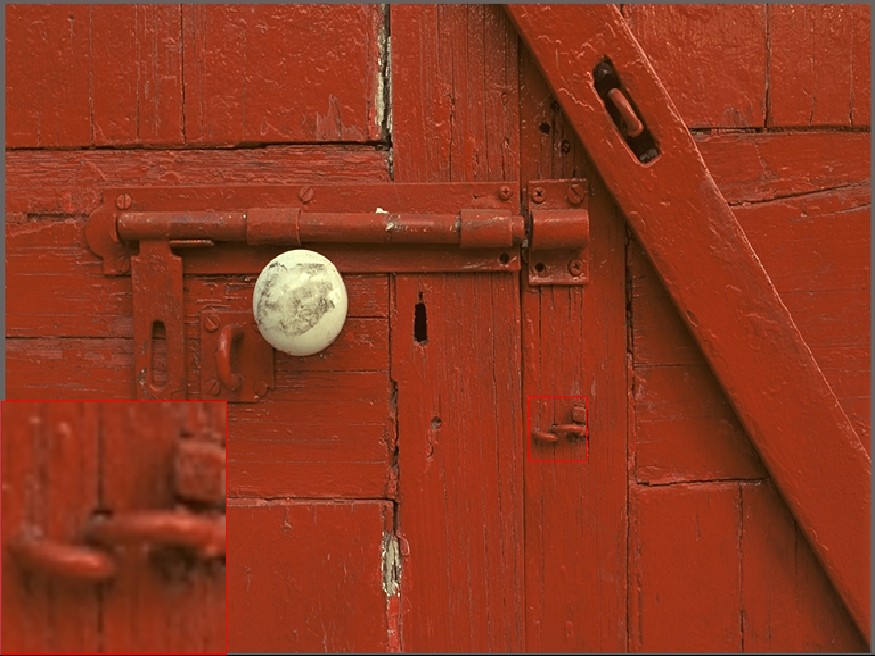}&	\includegraphics[width=0.13\linewidth]{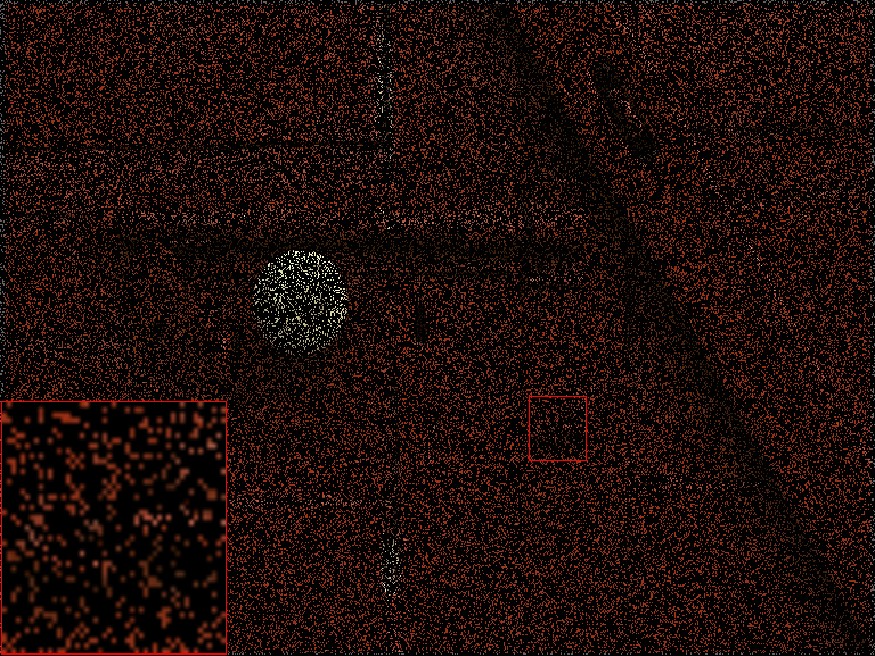}&	\includegraphics[width=0.13\linewidth]{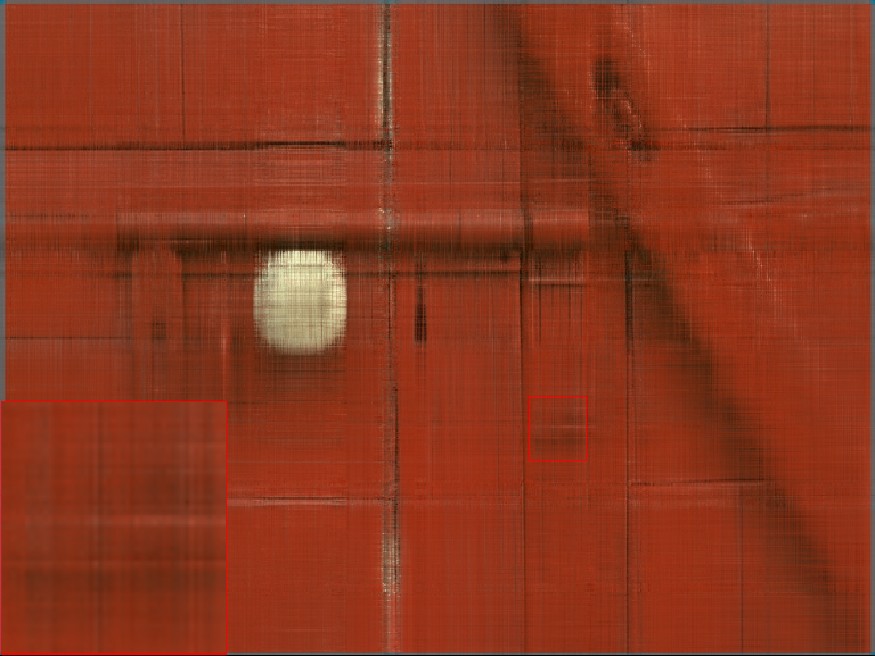}& 	\includegraphics[width=0.13\linewidth]{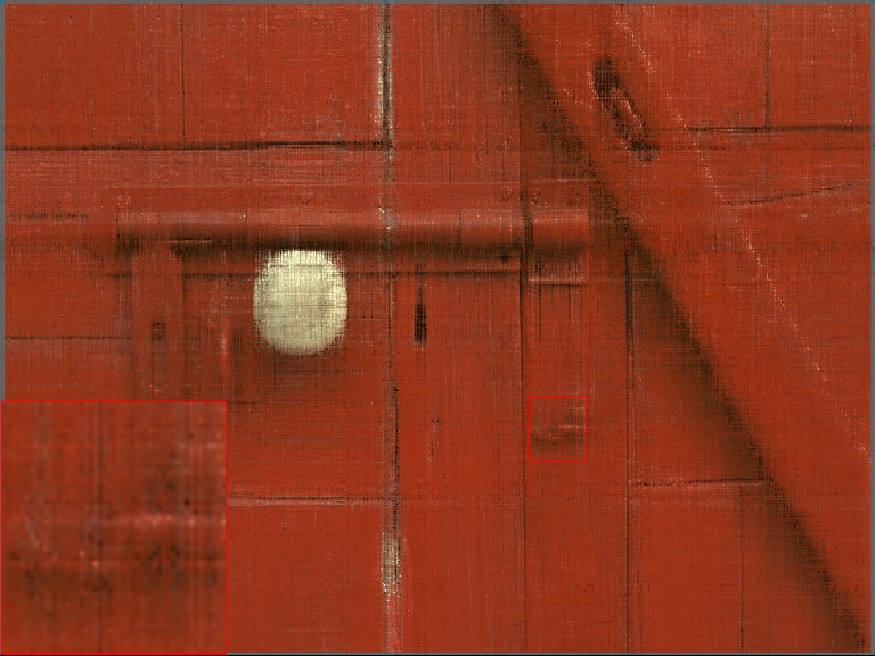}&	\includegraphics[width=0.13\linewidth]{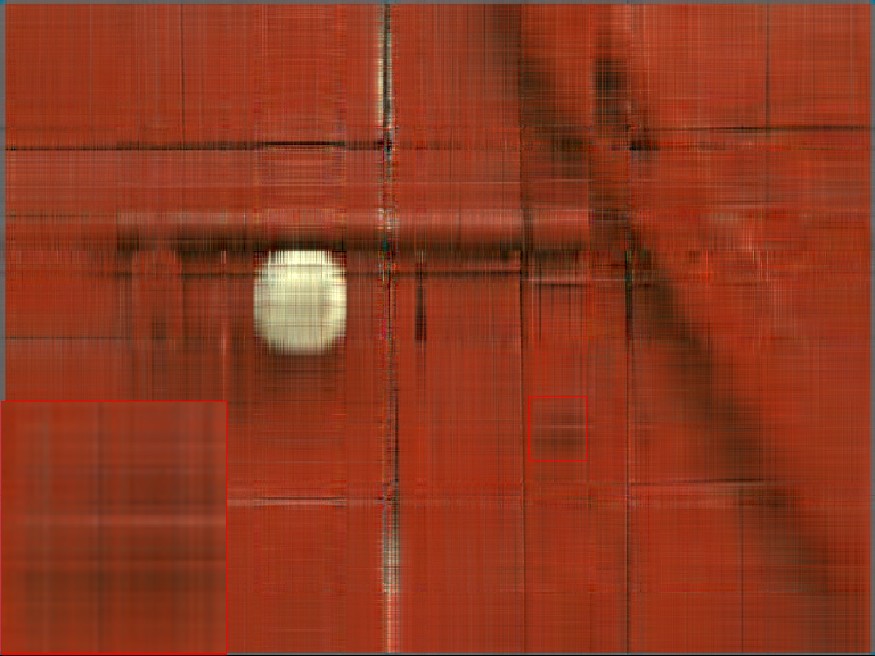}&\includegraphics[width=0.13\linewidth]{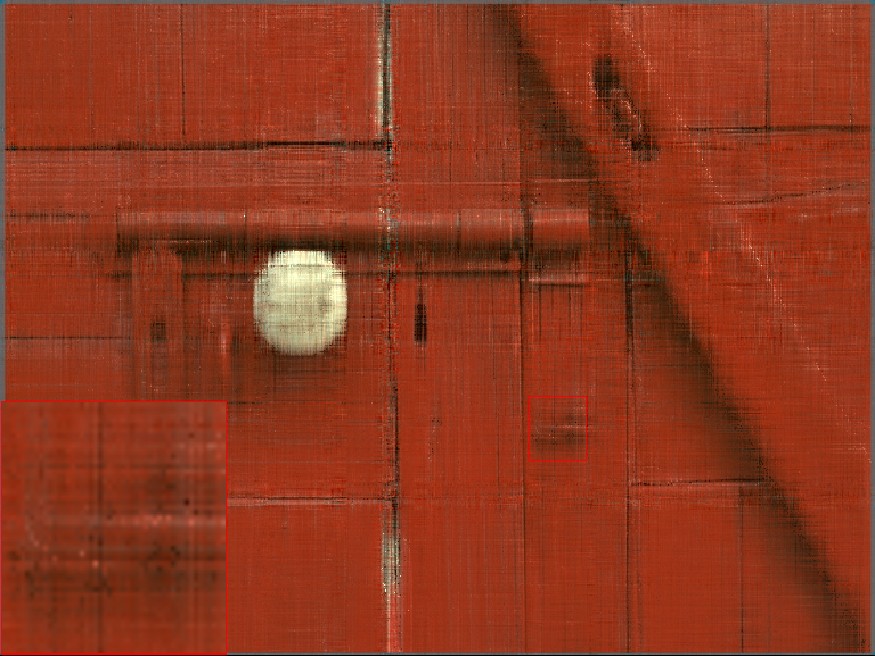}&		\includegraphics[width=0.13\linewidth]{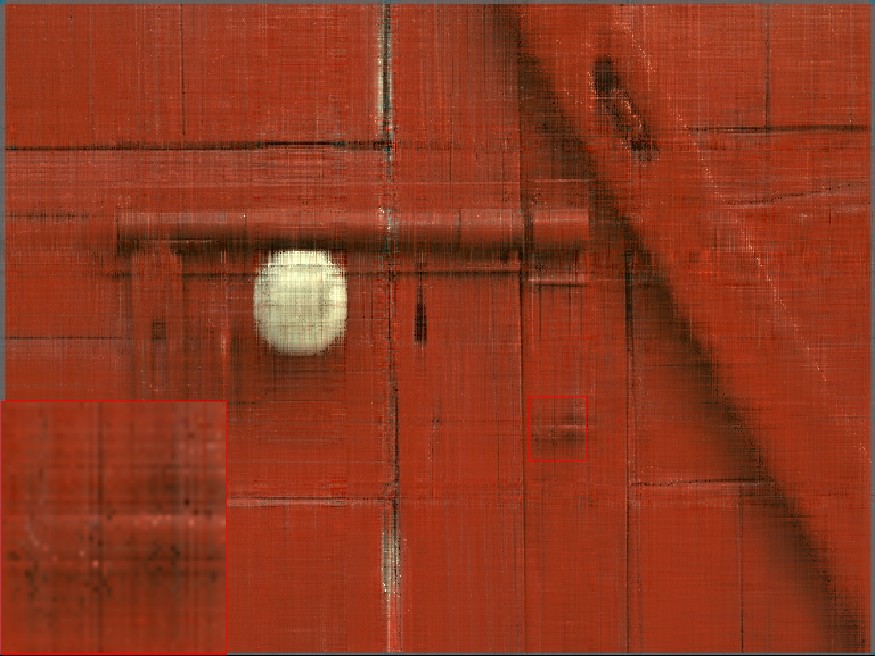}\\	\includegraphics[width=0.13\textwidth]{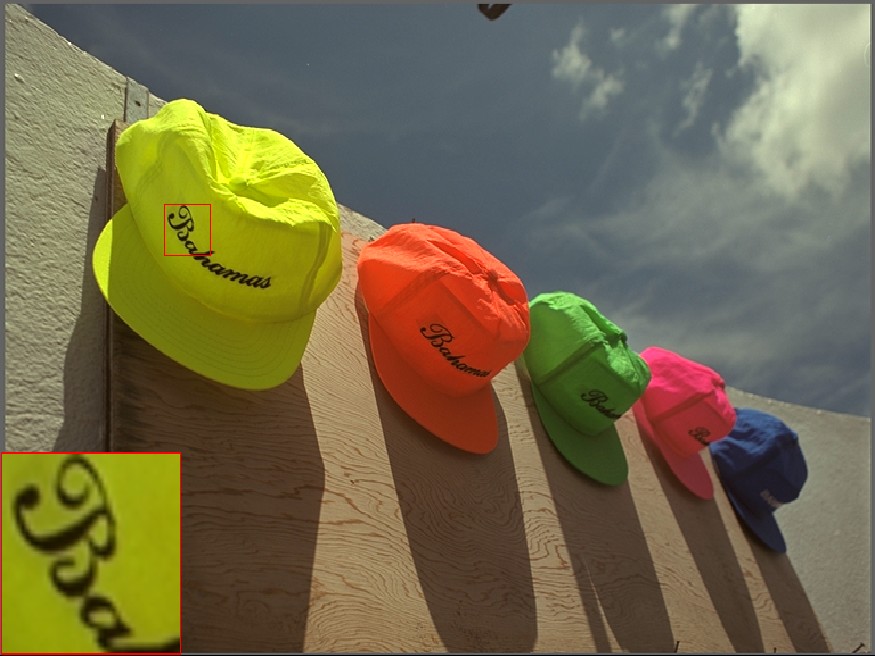}&	\includegraphics[width=0.13\linewidth]{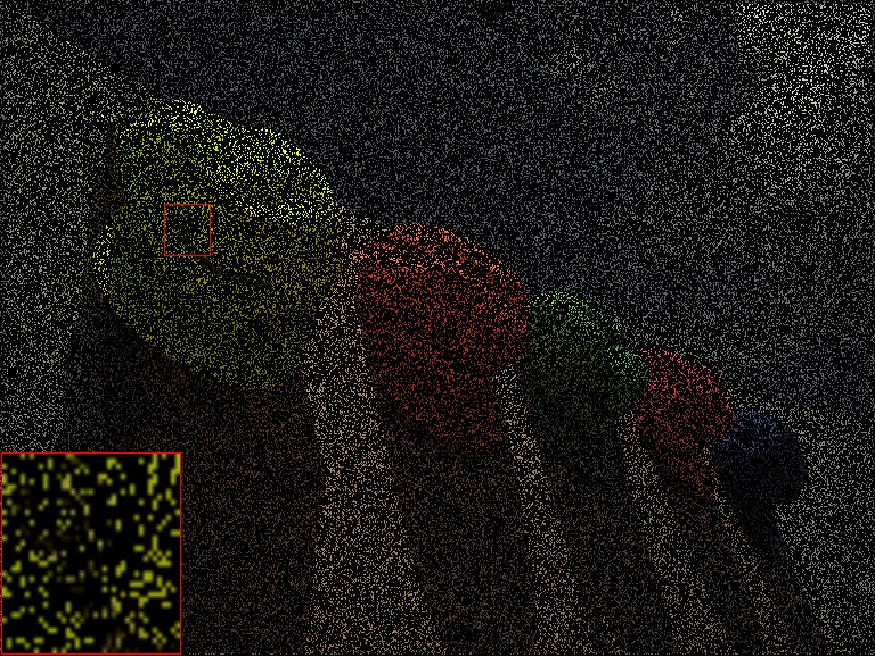}&  \includegraphics[width=0.13\linewidth]{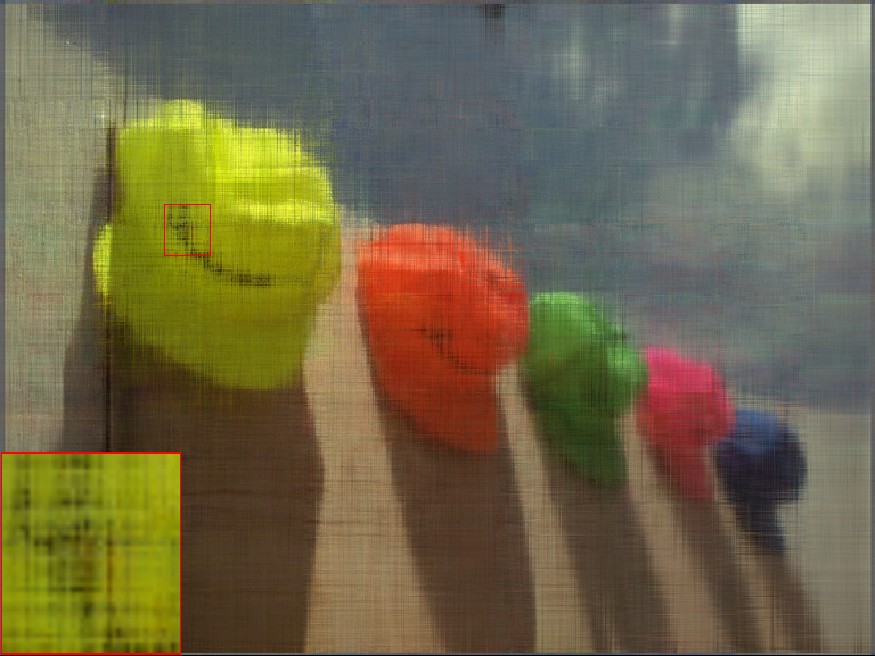}& 	\includegraphics[width=0.13\linewidth]{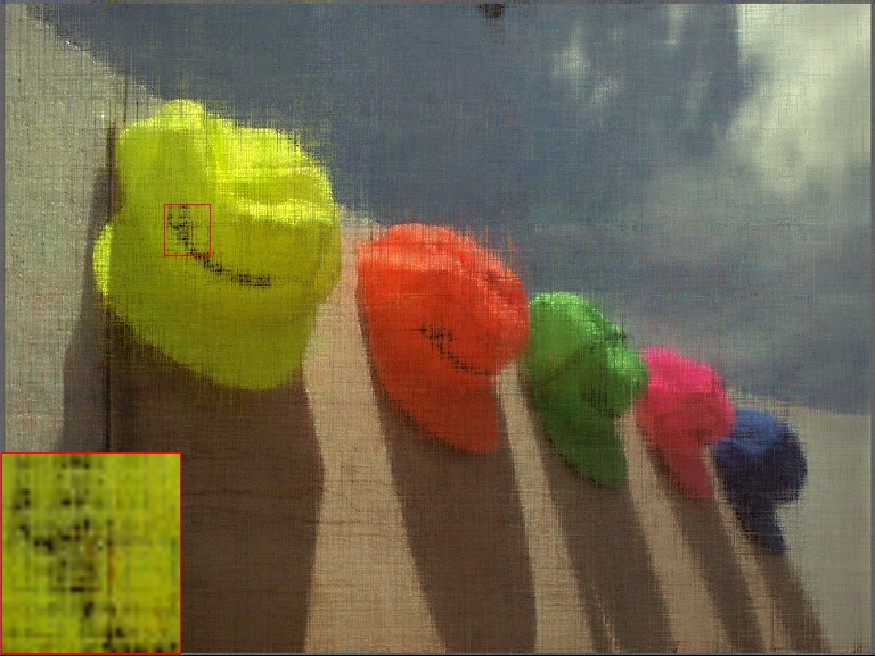}&	\includegraphics[width=0.13\linewidth]{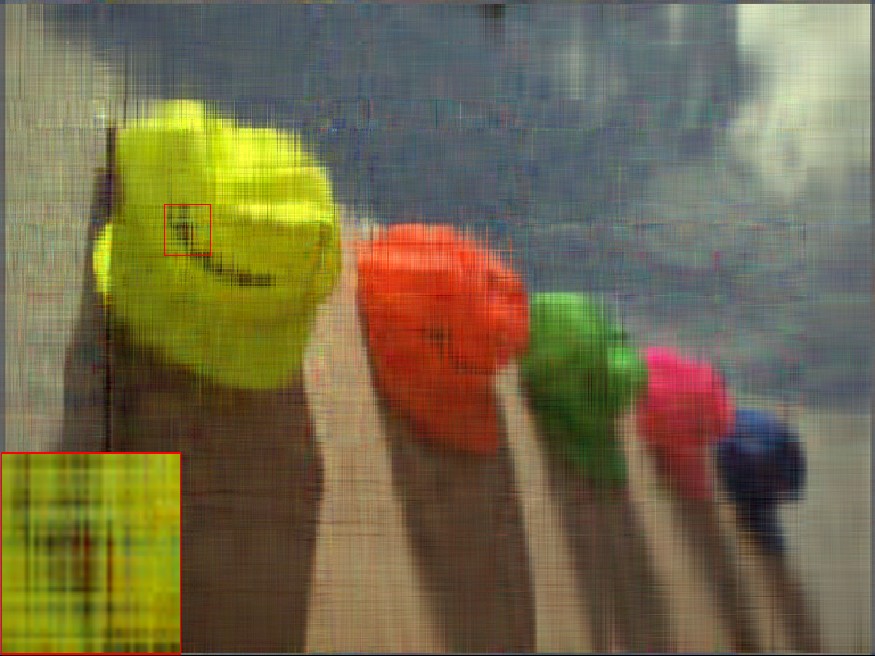}& \includegraphics[width=0.13\linewidth]{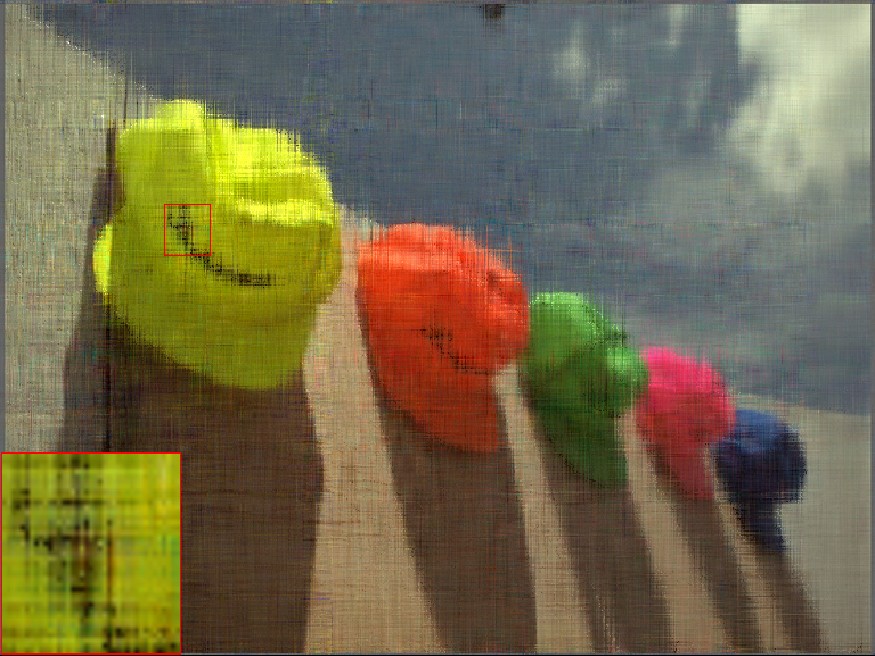}&	\includegraphics[width=0.13\linewidth]{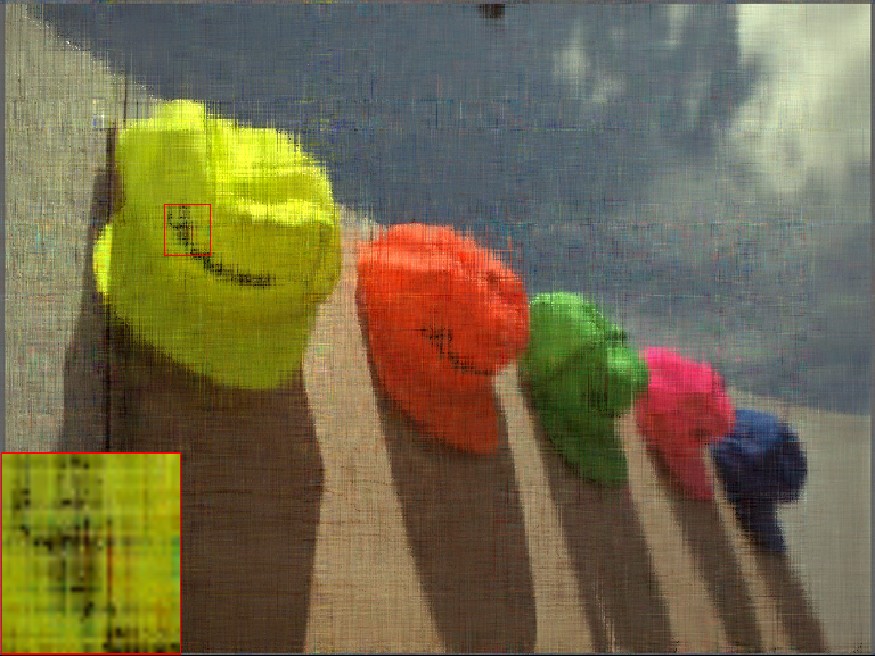}\\	\includegraphics[width=0.13\textwidth]{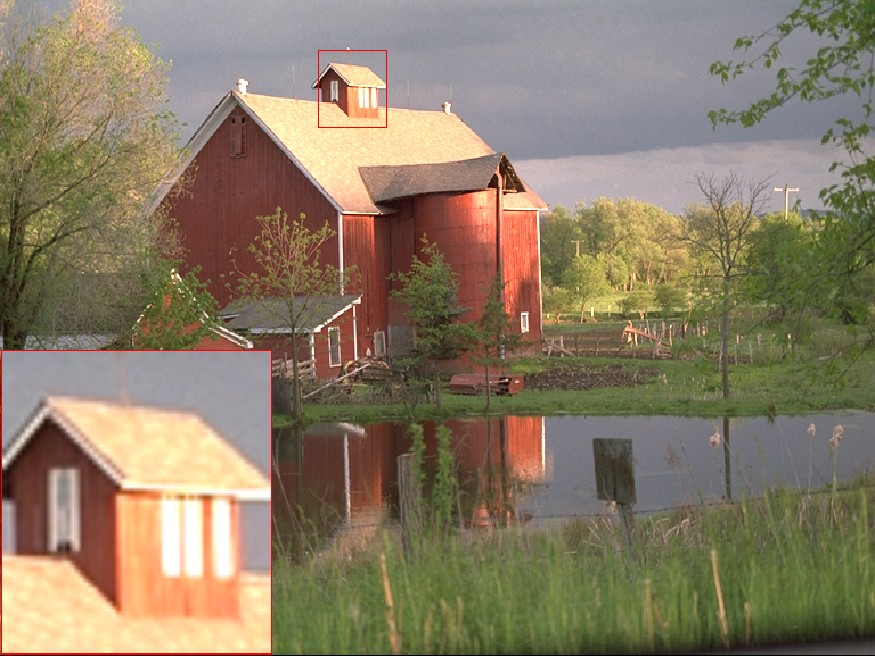}&	\includegraphics[width=0.13\linewidth]{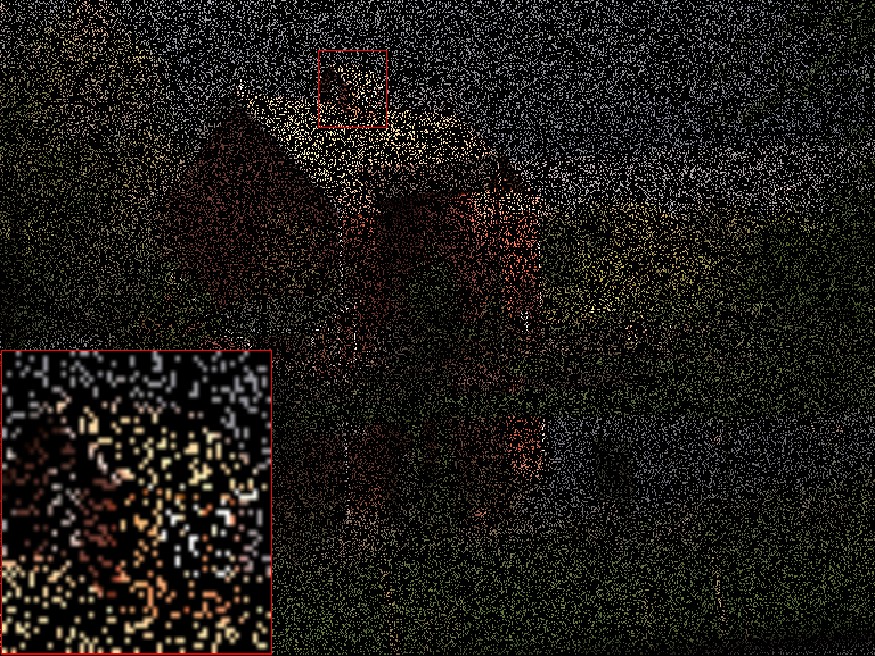}&  \includegraphics[width=0.13\linewidth]{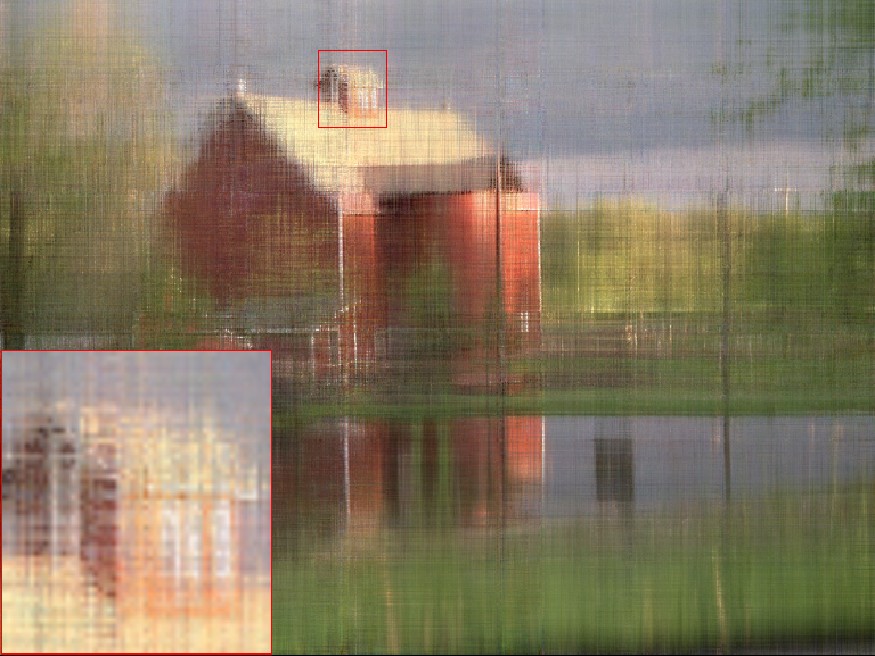}& 	\includegraphics[width=0.13\linewidth]{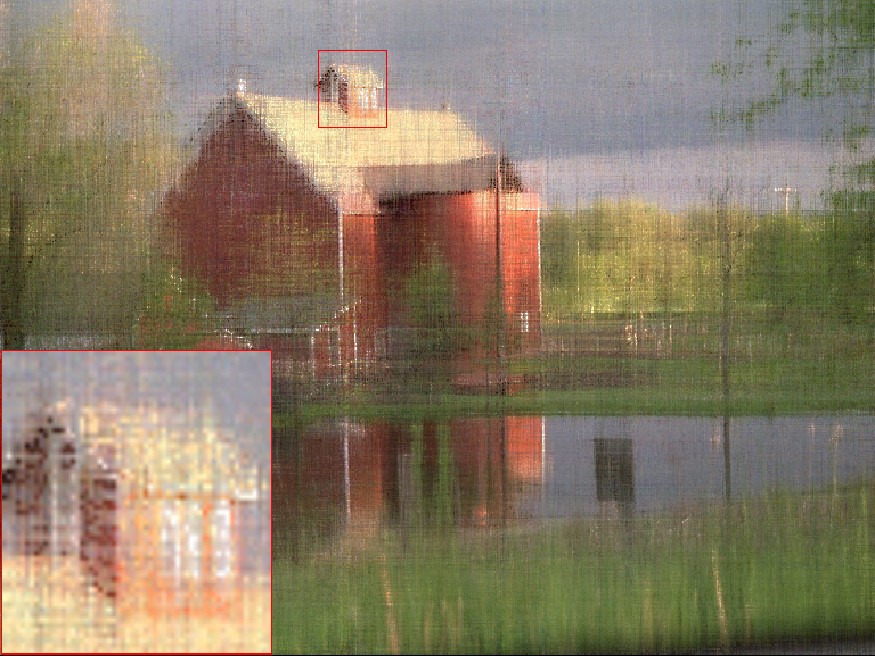}&	\includegraphics[width=0.13\linewidth]{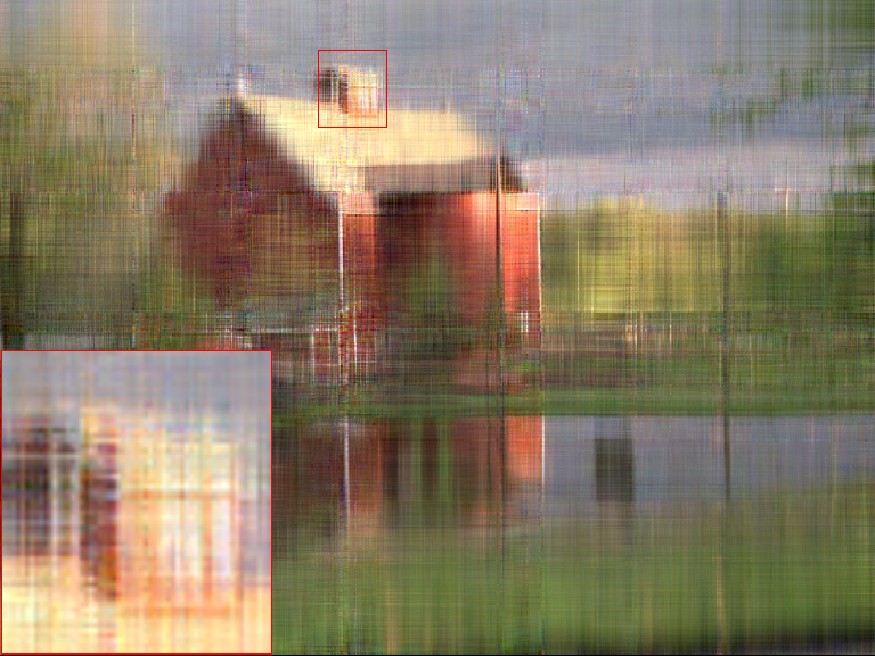}& \includegraphics[width=0.13\linewidth]{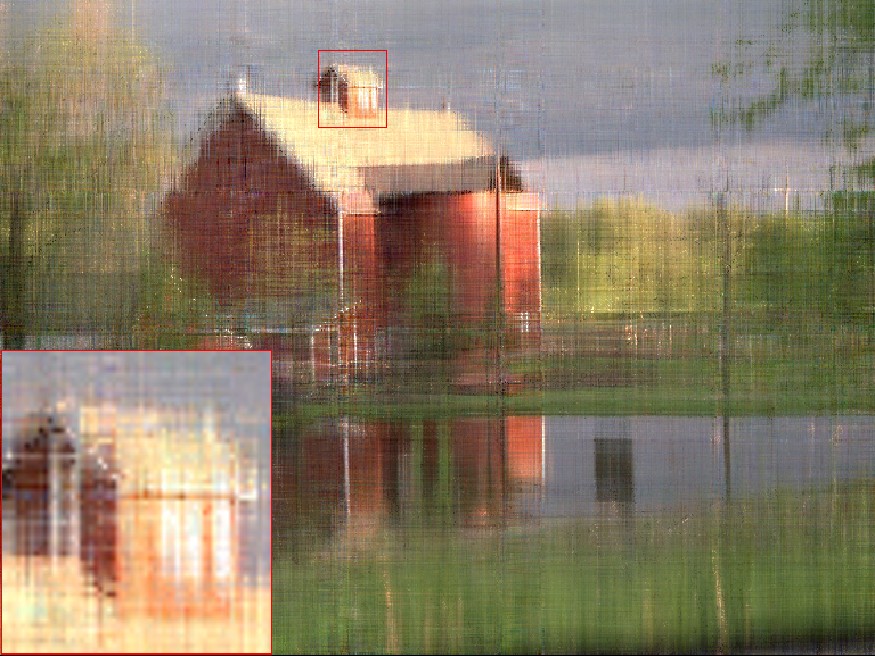}&	\includegraphics[width=0.13\linewidth]{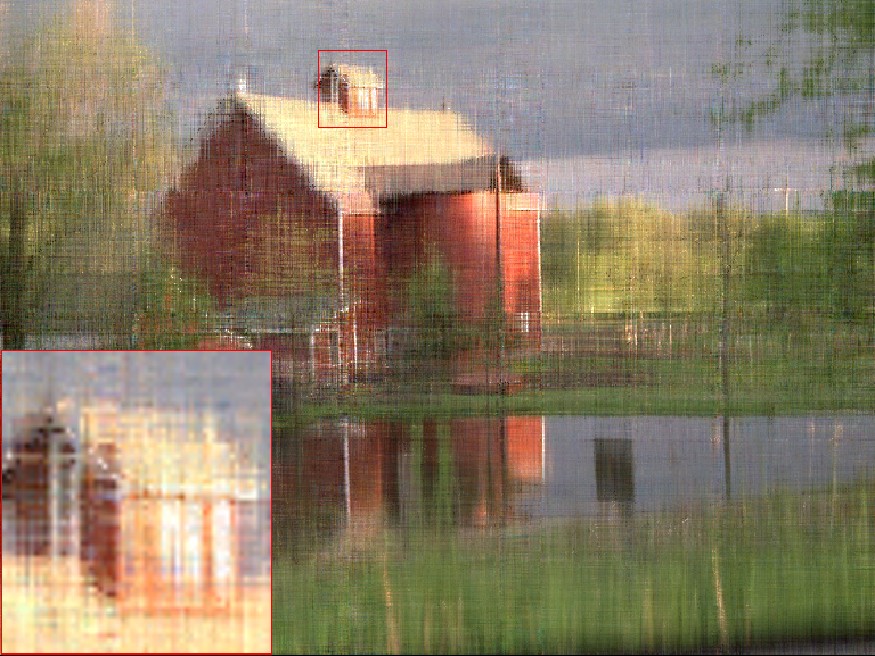}\\	\includegraphics[width=0.13\textwidth]{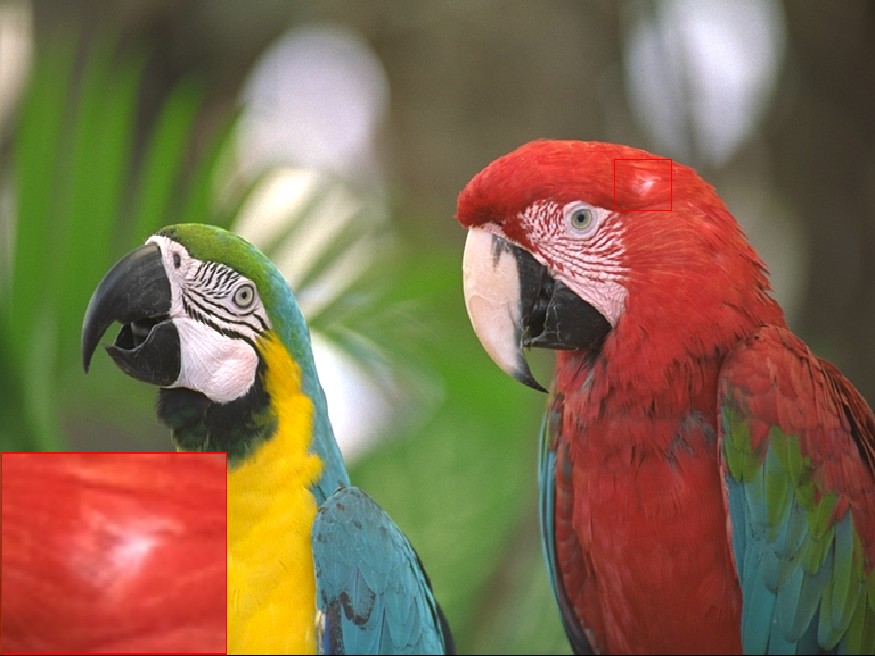}&	\includegraphics[width=0.13\linewidth]{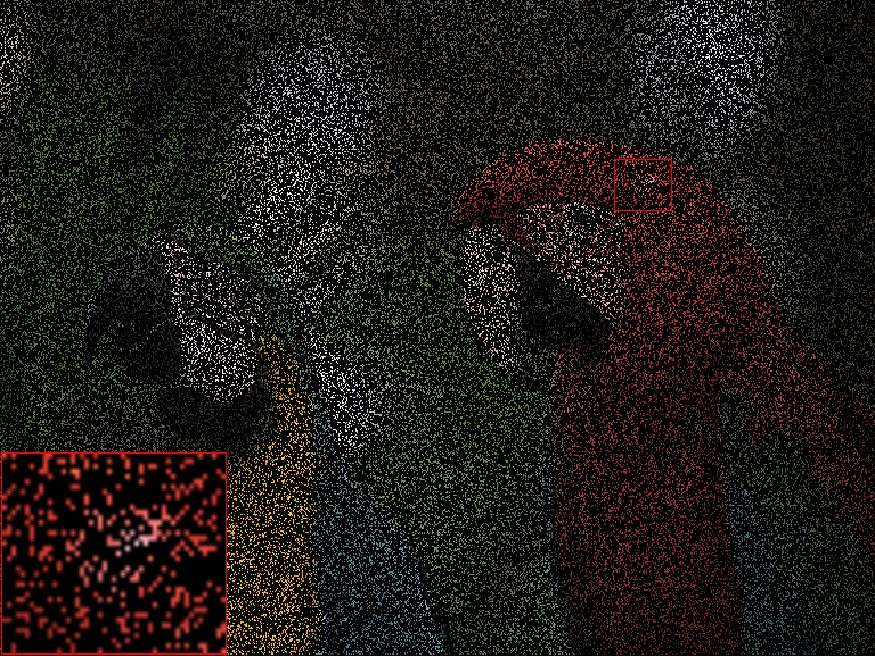}&  \includegraphics[width=0.13\linewidth]{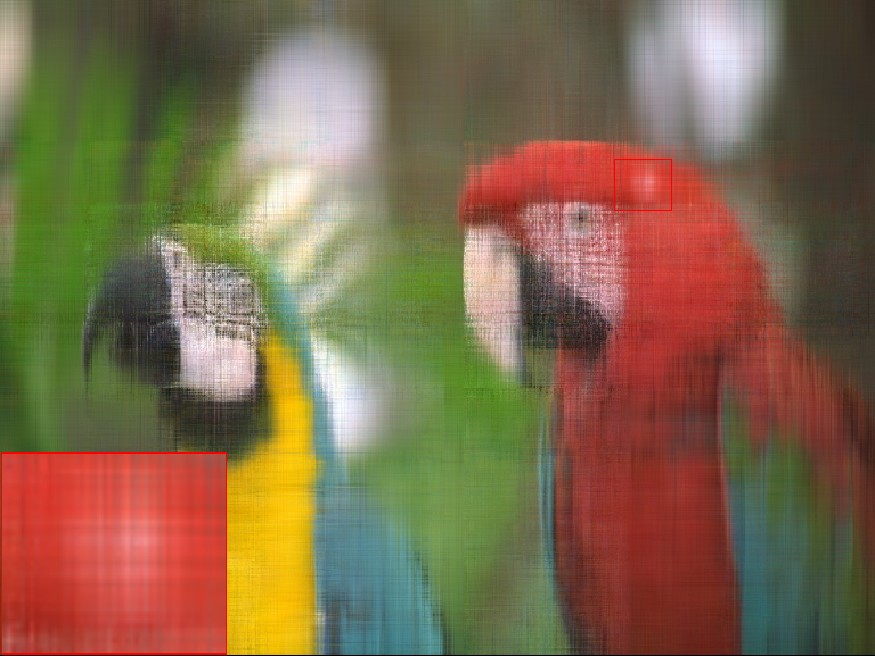}& 	\includegraphics[width=0.13\linewidth]{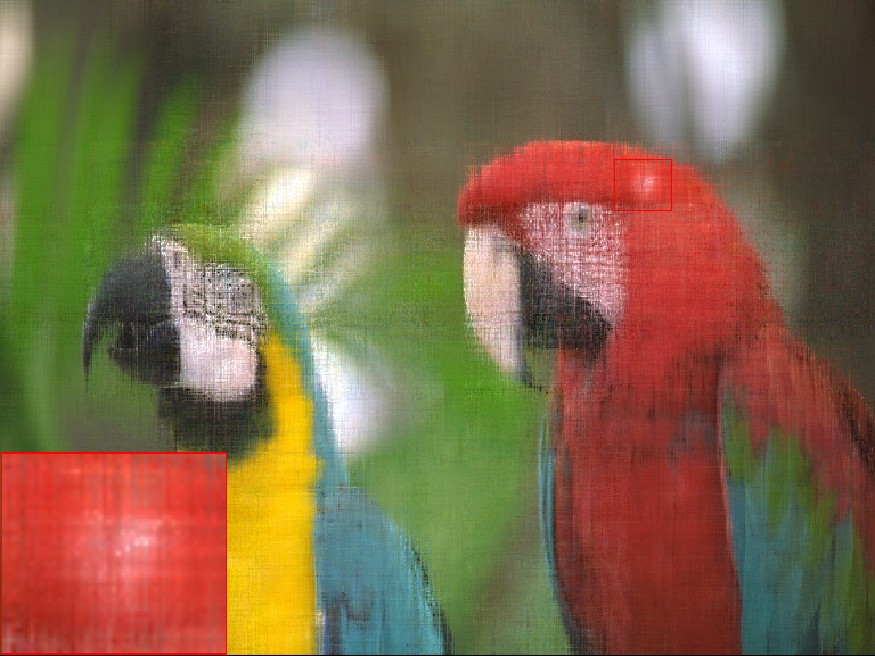}&	\includegraphics[width=0.13\linewidth]{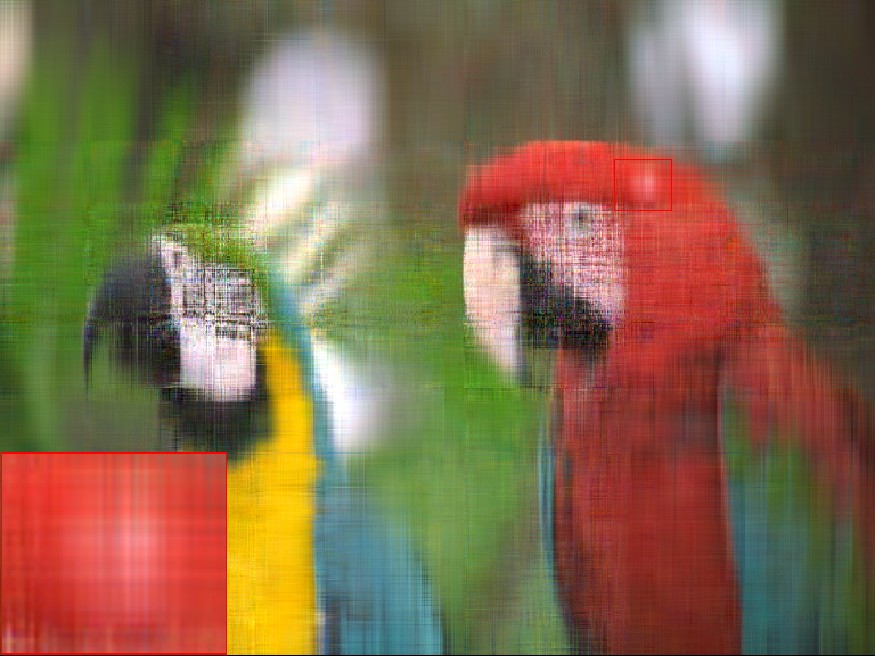}& \includegraphics[width=0.13\linewidth]{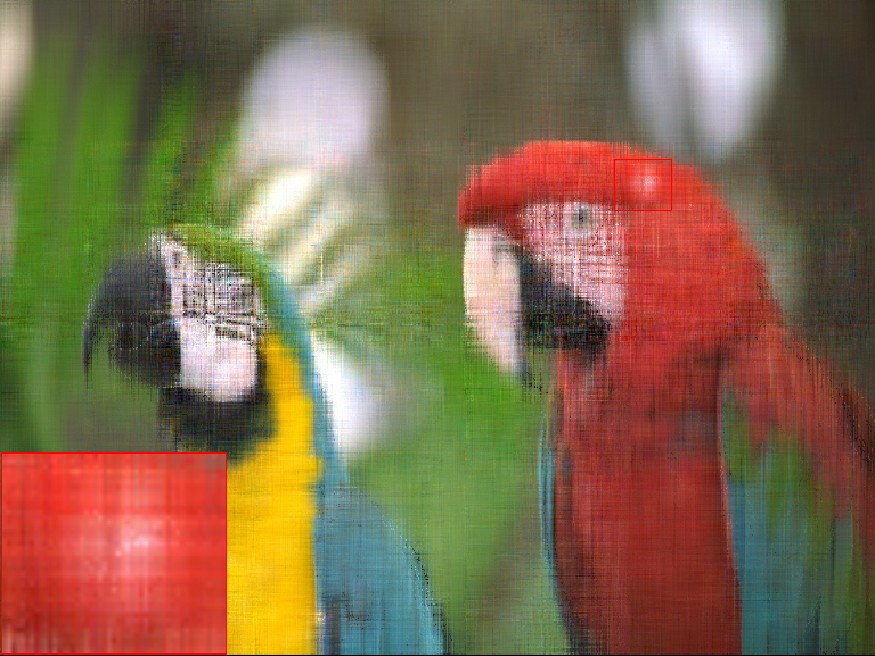}&	\includegraphics[width=0.13\linewidth]{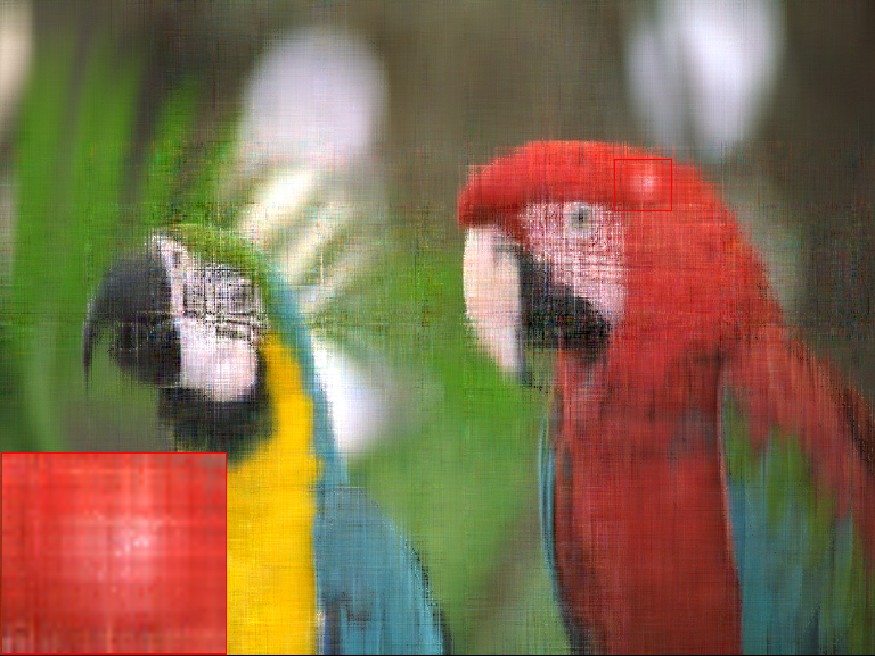}\\	\renewcommand\arraystretch{1}		Original & Observed &LRQA-1& DFN &QLNF&QMCUR$\_$length&QMCUR$\_$uniform	\end{tabular}	\caption{Visual comparison of various methods for color image recovery. From top to bottom: Image03, Image04, Image05, and Image06, respectively.} \label{tablecolorshow} 	\end{center}  \end{figure*}
From these results, the following observations are made. First, the PSNR results of the proposed methods (QMCUR\_length and QMCUR\_uniform) and three other methods (LRQA-1, QLNF, and Q-DFN) indicate that QMCUR\_length and QMCUR\_uniform achieve the best recovery performance in most cases. Furthermore, the proposed QMCUR\_uniform is the fastest, demonstrating its scalability on large-scale color images. As shown in Fig. \ref{tablecolorshow}, the visual effect of the QMCUR\_length and QMCUR\_uniform methods is superior to the comparative methods with MR = 80$\%$, which is consistent with our quantitative results. Both of these findings further demonstrate the superiority of our methods over comparative methods for handling large-scale data.
\section{Conclusion}

This work presented an efficient quaternion matrix CUR method for computing low-rank approximation of a quaternion data matrix. The QMCUR approximation offers a balance between accuracy and computational costs by selecting specific column and row submatrices of a given quaternion matrix. Next, we conducted a perturbation analysis of this approximation and concluded that the error in the quaternion spectral norm is correlated with the noise quaternion matrix in the first order. Computational experiments illustrated that the QMCUR approximation methods are significantly faster than comparative low-rank quaternion matrix approximation methods, without sacrificing the quality of reconstruction on both synthetic and color image datasets.

\bmhead{Acknowledgments}
	
This research is supported by University of Macau (MYRG2022-00108-FST, MYRG-CRG2022-00010-ICMS), the Science and Technology Development Fund, Macau S.A.R (0036/2021/AGJ), the National Key Research Project (2022YFB3904104), and the Key Project from National Natural Science Foundation of China (42230406).
\backmatter

\section*{Declarations}
\bmhead{Ethical Approval}
Not applicable.
\bmhead{Availability of supporting data}
The data sets generated during and/or analyzed during the current study are available
from the corresponding author on reasonable request.
\bmhead{Competing interests}
The authors declare no competing interests.
\bmhead{Funding}
This research is supported by University of Macau (MYRG2022-00108-FST, MYRG-CRG2022-00010-ICMS), the Science and Technology Development Fund, Macau S.A.R (0036/2021/AGJ), the National Key Research Project (2022YFB3904104), and the Key Project from National Natural Science Foundation of China (42230406).
\bmhead{Authors’ contributions}
Pengling Wu: Methodology, Software, Writing - Original Draft. 
Kit Ian Kou: Conceptualization, Validation, Resources. 
Hongmin Cai: Methodology, Validation, Writing - Review and Editing.
Zhaoyuan Yu: Formal analysis, Writing - Review and Editing.

\bibliography{Paper_Qrpca_cur}

\end{document}